\documentclass[final,onefignum,onetabnum]{siamart220329}

\usepackage{amsmath,amssymb,amsfonts}

\usepackage{amsmath} % for \dfrac
\usepackage{bm}
\usepackage{physics}
\usepackage{tikz,pgfplots}
\usetikzlibrary{angles,quotes} % for pic (angle labels)
\usetikzlibrary{calc}
\usetikzlibrary{decorations.markings}
\tikzset{>=latex} % for LaTeX arrow head

\usepackage{xcolor}

\usepackage{braket}

\usepackage{array}

\newtheorem{defn}{Definition}[section]
\newtheorem{thm}{Theorem}[section]
\newtheorem{remark}{Remark}[section]
\newtheorem{assumption}{Assumption}[section]

\newtheorem{prop}{Proposition}[section]
\renewcommand{\d}{\mathrm{d}}
\renewcommand{\L}{\mathrm{L}}
\renewcommand{\H}{\mathrm{H}}

\usepackage{algorithmic}

\usepackage{graphicx,epstopdf}

\Crefname{ALC@unique}{Line}{Lines}

\usepackage{amsopn}

\pgfplotsset{compat=1.18}

\usepackage{xspace}
\usepackage{bold-extra}
\usepackage[most]{tcolorbox}

\colorlet{texcscolor}{blue!50!black}
\colorlet{texemcolor}{red!70!black}
\colorlet{texpreamble}{red!70!black}
\colorlet{codebackground}{black!25!white!25}

\lstdefinestyle{siamlatex}{%
  style=tcblatex,
  texcsstyle=*\color{texcscolor},
  texcsstyle=[2]\color{texemcolor},
  keywordstyle=[2]\color{texemcolor},
  moretexcs={cref,Cref,maketitle,mathcal,text,headers,email,url},
}

\tcbset{%
  colframe=black!75!white!75,
  coltitle=white,
  colback=codebackground, % bottom/left side
  colbacklower=white, % top/right side
  fonttitle=\bfseries,
  arc=0pt,outer arc=0pt,
  top=1pt,bottom=1pt,left=1mm,right=1mm,middle=1mm,boxsep=1mm,
  leftrule=0.3mm,rightrule=0.3mm,toprule=0.3mm,bottomrule=0.3mm,
  listing options={style=siamlatex}
}

\newtcbinputlisting[use counter=example]{\examplefile}[3][]{%
  title={Example~\thetcbcounter: #2},listing file={#3},#1}

\DeclareTotalTCBox{\code}{ v O{} }
{ %fontupper=\ttfamily\color{texemcolor},
  fontupper=\ttfamily\color{black},
  nobeforeafter,
  tcbox raise base,
  colback=codebackground,colframe=white,
  top=0pt,bottom=0pt,left=0mm,right=0mm,
  leftrule=0pt,rightrule=0pt,toprule=0mm,bottomrule=0mm,
  boxsep=0.5mm,
  #2}{#1}

\patchcmd\newpage{\vfil}{}{}{}
\begin{tcbverbatimwrite}{tmp_\jobname_header.tex}
\title{Linear-Quadratic optimal control for boundary controlled networks of waves\thanks{Submitted to the editors on February 21, 2024.
\funding{This research was conducted with the financial supports of the F.R.S-FNRS (Belgium) and the German Research Foundation (DFG). Anthony Hastir was supported by a FNRS Postdoctoral Fellowship, Grant CR 40010909. He has also been supported by a FNRS mobility funding, reference 40020129. He is now supported by the DFG via the Grant HA 10262/2-1.}}
}

\author{Anthony Hastir\thanks{University of Wuppertal, School of Mathematics and Natural Sciences
Gaußstraße 20, 42119 Wuppertal, Germany, \email{hastir@uni-wuppertal.de}.}
\and Birgit Jacob\thanks{University of Wuppertal, School of Mathematics and Natural Sciences
Gaußstraße 20, 42119 Wuppertal, Germany, \email{bjacob@uni-wuppertal.de}.}
\and Hans Zwart\thanks{Department of Applied Mathematics, University of Twente, P.O. Box 217, 7500 AE, Enschede, The Netherlands and Department of Mechanical Engineering, Eindhoven University of Technology, P.O. Box 513, 5600 MB, Eindhoven, The Netherlands, \email{h.j.zwart@utwente.nl}}
}

\headers{LQ optimal control for BC networks of waves}{Anthony Hastir, Birgit Jacob, and Hans Zwart}
\end{tcbverbatimwrite}
\title{Linear-Quadratic optimal control for boundary controlled networks of waves\thanks{Submitted to the editors on February 21, 2024.
\funding{This research was conducted with the financial supports of the F.R.S-FNRS (Belgium) and the German Research Foundation (DFG). Anthony Hastir was supported by a FNRS Postdoctoral Fellowship, Grant CR 40010909. He has also been supported by a FNRS mobility funding, reference 40020129. He is now supported by the DFG via the Grant HA 10262/2-1.}}
}

\author{Anthony Hastir\thanks{University of Wuppertal, School of Mathematics and Natural Sciences
Gaußstraße 20, 42119 Wuppertal, Germany, \email{hastir@uni-wuppertal.de}.}
\and Birgit Jacob\thanks{University of Wuppertal, School of Mathematics and Natural Sciences
Gaußstraße 20, 42119 Wuppertal, Germany, \email{bjacob@uni-wuppertal.de}.}
\and Hans Zwart\thanks{Department of Applied Mathematics, University of Twente, P.O. Box 217, 7500 AE, Enschede, The Netherlands and Department of Mechanical Engineering, Eindhoven University of Technology, P.O. Box 513, 5600 MB, Eindhoven, The Netherlands, \email{h.j.zwart@utwente.nl}}
}

\headers{LQ optimal control for BC networks of waves}{Anthony Hastir, Birgit Jacob, and Hans Zwart}

\ifpdf
\hypersetup{ pdftitle={Linear-Quadratic optimal control for boundary controlled boundary observed hyperbolic PDEs} }
\fi

\begin{document}
\maketitle

%% ------------------------------------------------------------------
%% ABSTRACT
%% ------------------------------------------------------------------
\begin{tcbverbatimwrite}{tmp_\jobname_abstract.tex}
\begin{abstract}
Linear-Quadratic optimal controls are computed for a class of boundary controlled, boundary observed hyperbolic infinite-dimensional systems, which may be viewed as networks of waves. The main results of this manuscript consist in converting the infinite-dimensional continuous-time systems into infinite-dimensional discrete-time systems for which the operators dynamics are matrices, in solving the LQ-optimal control problem in discrete-time and then in interpreting the solution in the continuous-time variables, giving rise to the optimal boundary control input. The results are applied to two examples, a small network of three vibrating strings and a co-current heat-exchanger, for which boundary sensors and actuators are considered. 
\end{abstract}


\end{tcbverbatimwrite}
\begin{abstract}
Linear-Quadratic optimal controls are computed for a class of boundary controlled, boundary observed hyperbolic infinite-dimensional systems, which may be viewed as networks of waves. The main results of this manuscript consist in converting the infinite-dimensional continuous-time systems into infinite-dimensional discrete-time systems for which the operators dynamics are matrices, in solving the LQ-optimal control problem in discrete-time and then in interpreting the solution in the continuous-time variables, giving rise to the optimal boundary control input. The results are applied to two examples, a small network of three vibrating strings and a co-current heat-exchanger, for which boundary sensors and actuators are considered.
\end{abstract}



%% ------------------------------------------------------------------
%% END HEADER
%% ------------------------------------------------------------------

\section{Introduction}
\label{sec:intro}
The Linear-Quadratic (LQ) optimal control problem is a well-known and well-established problem in control theory. For continuous-time finite-dimensional systems driven by the following dynamics
\begin{align}
  \dot{x}(t) &= Ax(t) + Bu(t),\label{State_Finite}\\
  x(0) &= x_0,\label{Init_Finite}\\
  y(t) &= Cx(t) + Du(t),\label{Output_Finite}
\end{align}
with $A\in\mathbb{C}^{n\times n}, B\in\mathbb{C}^{p\times n}, C\in\mathbb{C}^{m\times n}, D\in\mathbb{C}^{m\times p}$, the optimal control input that minimizes the cost functional
\begin{equation}
J(x_0,u) = \int_0^\infty \left[\Vert u(t)\Vert^2 + \Vert y(t)\Vert^2\right]\d t
\label{Cost}
\end{equation}
is given by the state feedback $u(t) = Kx(t) = -B^*\Theta x(t)$, where $\Theta = \Theta^*$ is the smallest nonnegative solution of the control algebraic Riccati equation
\begin{equation}
A^*\Theta + \Theta A + C^*C = \left(\Theta B + C^*D\right)\left(I + D^*D\right)^{-1}\left(B^*\Theta + D^*C\right).
\label{Riccati_Finite}
\end{equation}
This is valid provided that the so-called \textit{finite-cost condition} is satisfied, that is, for every initial condition $x_0$, there exists an input $u$ such that the cost \eqref{Cost} is finite. When $A$ is an (unbounded) operator on a Hilbert space $X$ that generates a strongly continuous semigroup and $B, C$ and $D$ are bounded linear operators, that is, $B\in\mathcal{L}(U,X), C\in\mathcal{L}(X,Y), D\in\mathcal{L}(U,Y)$ with the Hilbert spaces $U$ and $Y$ being the input and the output spaces, respectively, the solution to the optimal control problem is described in a similar way as in the finite-dimensional case. If the finite-cost condition is satisfied, the optimal solution is still given by the state feedback $u(t) = -B^*\Theta x(t)$ with $\Theta\in\mathcal{L}(X)$ being the smallest nonnegative solution to the operator Riccati equation, i.e., for all $x\in D(A)$, $\Theta$ satisfies
\begin{equation}
\langle Ax, \Theta x\rangle + \langle \Theta x, Ax\rangle + \langle Cx, Cx\rangle = \langle (I + D^*D)^{-1}(B^*\Theta + D^*C)x,(B^*\Theta + D^*C)x\rangle,
\label{Riccati_WeakForm}
\end{equation}
see \cite{Curtain_Pritchard_Riccati}, \cite[Chapter 4]{CurtainPritchard}, \cite{Kwakernaak_Sivan}. In both cases, the closed-loop system is obtained by connecting the optimal control input to the original system \eqref{State_Finite}--\eqref{Output_Finite}. Moreover, the optimal cost satisfies
\begin{equation*}
  \langle \Theta x_0,x_0\rangle = \inf_{u\in\L^2(0,\infty;U)}J(x_0,u).
\end{equation*}
These results break down in their original forms when the input and output operators are allowed to be unbounded. As is pointed out in \cite{Weiss_Zwart}, the Riccati equation \eqref{Riccati_WeakForm} has to be adapted mostly because $\Theta x\notin D(B^*)$ for $x\in D(A)$. For that reason, extensions of the operator $B^*$ have to be considered. The authors in \cite{WeissWeiss} and \cite{Staffans_LQ_1998} proposed the following operator Riccati equation
\begin{equation}
  \langle Ax, \Theta x\rangle + \langle \Theta x, Ax\rangle + \langle Cx, Cx\rangle = \langle (\Omega^*\Omega)^{-1}(B^*_w\Theta + D^*C)x,(B^*_w\Theta + D^*C)x\rangle,
\label{Riccati_Unbounded}
\end{equation}
$x\in D(A)$, with $B^*_w$ being the Yosida extension of the operator $B^*$, see e.g.~\cite{Weiss_Yosida_Ext} or \cite{Tucsnak_Weiss_LTI_Case} for the definition of that operator. Another point of difference between \eqref{Riccati_WeakForm} and \eqref{Riccati_Unbounded} is the operator $\Omega^*\Omega$, which can be different from $I + D^*D$ because of the unbounded nature of $B$ and $C$, see e.g.~\cite{Opmeer_MTNS} for a counterexample highlighting this feature. In addition, the closed-loop system is not necessarily obtained by connecting the optimal feedback to the original system. Those problems are emphasized in \cite{Weiss_Zwart} in which a counterexample is considered. Different variations of the Riccati equation \eqref{Riccati_Unbounded} exist, see e.g.~\cite{Pritchard_Salamon}, \cite{WeissWeiss}, \cite{Curtain_LQ_SIAM} and \cite{mikkola2016riccati}. In \cite{Curtain_LQ_SIAM}, Riccati equations for well-posed linear systems have also been considered, under the assumption that $0$ is in the resolvent set of the operator dynamics, giving the possibility to define a reciprocal Riccati equation, based on the reciprocal system associated to the original one. A few years later, \cite{mikkola2016riccati} proposed generalizations of previous works on LQ-optimal control by introducing an integral Riccati equation, able to treat a large class of infinite-dimensional dynamical systems with a large class of quadratic cost functions. More recently, another type of Riccati equation called the \textit{operator node Riccati equation} has been considered in \cite{Opmeer_MTNS} and \cite{Opmeer_Staffans} for general infinite-dimensional systems. Moreover, links between the solution of the LQ-optimal control problem for general PDEs and for the discrete-time system obtained by applying the internal Cayley transform to them are presented in \cite{Opmeer_Staffans}. However, obtaining a closed form for the solution to those Riccati equations, and in particular an operator $\Theta$ solution of \eqref{Riccati_Unbounded}, is challenging and almost impossible, except for very special cases. This is mainly because \eqref{Riccati_Unbounded} is a quadratic operator equation where all the operators are unbounded. This makes the results presented in \cite{WeissWeiss}, \cite{Opmeer_MTNS}, \cite{Opmeer_Staffans} very technical and difficult to implement in applications. In addition to the aforementioned results based on state-space descriptions, a frequency domain approach of LQ-optimal control exists and is known as the spectral factorization method, which consists in finding the standard spectral factor of a power spectral density called the Popov function. This method is extensively developed for infinite-dimensional systems with bounded input and output operators and finite-dimensional input and output spaces in \cite{CallierWinkin1990} and \cite{CallierWinkin1992}. In \cite{WeissWeiss}, the spectral factorization method is extended to unbounded control and observation operators, giving rise to a solution in the frequency domain. These results are applied to an example dealing with LQ-optimal control for a hyperbolic system with unbounded control operator in \cite{Weiss_Zwart}. However, finding spectral factors of the Popov function that satisfy appropriate regularity assumptions may be an extremely difficult task. Getting explicit/computable formulas for the spectral factors is as challenging as finding a solution to operator Riccati equations like \eqref{Riccati_Unbounded}.

In this manuscript, we solve the LQ-optimal control problem for a class of boundary controlled and boundary observed hyperbolic partial differential equations (PDEs) that may be viewed as networks of waves. In particular, our main result provides an explicit formula for the optimal control. Thanks to an equivalent representation of the system as an infinite-dimensional discrete-time system, our approach to treat LQ-optimal control bypasses the computation of a solution to the operator Riccati equation \eqref{Riccati_Unbounded}. This has the great advantage of providing an easy way to express the solution to the LQ-optimal control problem. Moreover, as the operators in the discrete-time systems describing the original PDEs are matrices, the Riccati equation associated to the LQ-optimal control problem is a matrix Riccati equation, easily solvable with mathematical softwares.

Hyperbolic infinite-dimensional systems constitute an important class of dynamical systems. Their study in the context of semigroup theory may be found in the early references \cite{Hyperbolic_Phillips_1957} and \cite{Hyperbolic_Phillips_1959}. Many physical applications in engineering in the broad sense may be modeled by hyperbolic partial differential equations (PDEs). For instance, this is the case for the system that describes the propagation of current and voltage in transmission lines, modeled by the telegrapher equations. The propagation of water and sediments in open channels is described by hyperbolic PDEs as well, and is known as the Saint-Venant-Exner equations. In chemical engineering, the dynamics of plug-flow tubular reactors are also governed by this type of equations. Furthermore, hyperbolic PDEs are able to describe the motion of an inviscid ideal gas in a rigid cylindrical pipe, via the Euler equations. More generally, conservation laws may be written via hyperbolic PDEs. Beyond these few examples, a comprehensive overview of many applications governed by hyperbolic PDEs may be found in \cite{BastinCoron_Book}. Therein, these equations are presented both in their nonlinear and linear forms, where aspects like stability and control are described as well. 

LQ optimal control for hyperbolic PDEs has already attracted a lot of attention in the last few years. For instance, this control problem has been considered in \cite{Aksikas_TAC} for a plug-flow tubular reactor with bounded input and output operators. The LQ optimal control problem has then been solved via spectral factorization. Another reference combining LQ optimal control and hyperbolic PDEs is \cite{Aksikas_Automatica_PDE-ODE}. Therein, the authors consider a class of hyperbolic PDEs coupled to ODEs via the boundary. They use the Fattorini approach \cite{BCS_Fattorini} to convert the abstract system with unbounded input operator to a system in which the input operator becomes bounded. They solve the LQ optimal control problem after this transformation. Few years later, a similar problem has been considered in \cite{Dehaye} for a general class of infinite-dimensional systems with unbounded input and output operators. The Fattorini approach consisting in adding an extra state variable to replace the control operator by a bounded operator has been used to tackle LQ optimal control. Moreover, the authors consider a Yosida-type approximation of the output operator, which makes it bounded as well.

Here, we focus on the following class of hyperbolic PDEs, sometimes called networks of waves
\begin{align}
\frac{\partial \tilde{z}}{\partial t}(\zeta,t) &= -\frac{\partial}{\partial \zeta}(\lambda_0(\zeta)\tilde{z}(\zeta,t)) + M(\zeta)\tilde{z}(\zeta,t),\nonumber\\
\tilde{z}(\zeta,0) &= \tilde{z}_0(\zeta).\label{StateEquation_Diag_Uniform}
\end{align}
The variables $t\geq 0$ and $\zeta\in [0,1]$ are called time and space variables, respectively, whereas $\tilde{z}(\cdot,t)\in\L^2(0,1;\mathbb{C}^n) =: X$ is the state with initial condition $\tilde{z}_0\in X$. The PDEs are complemented with the following inputs and outputs
\begin{align}
\left[\begin{matrix}0\\ I\end{matrix}\right]u(t) &= -\lambda_0(0) K\tilde{z}(0,t) - \lambda_0(1) L \tilde{z}(1,t),\label{Input_Diag_Uniform}\\
y(t) &= -\lambda_0(0) K_y \tilde{z}(0,t) - \lambda_0(1) L_y \tilde{z}(1,t),\label{Output_Diag_Uniform}
\end{align}
$t\geq 0$, where $u(t)\in\mathbb{C}^p, p\in\mathbb{N}$ is the control input and $y(t)\in\mathbb{C}^m, m\in\mathbb{N}$ is the output. The scalar-valued function $\lambda_0\in \L^\infty(0,1;\mathbb{R}^+)$ satisfies $\lambda_0(\zeta)\geq \epsilon > 0$ a.e.~on $[0,1]$. $M(\zeta)$ is a space-dependent square matrix that satisfies $M\in\L^\infty(0,1;\mathbb{C}^{n\times n})$ while $K\in\mathbb{C}^{n\times n}, L\in\mathbb{C}^{n\times n}, K_y\in\mathbb{C}^{m\times n}$ and $L_y\in\mathbb{C}^{m\times n}$. 

This class of hyperbolic systems appears when looking at networks of conservation laws such as networks of electric lines, chains of density-velocity systems or genetic regulatory networks to mention a few, see e.g. \cite{BastinCoron_Book} for a comprehensive overview of the covered applications. We also refer the reader to the book \cite[Chapter 6]{Luo_Guo} in which a more general version of the class \eqref{StateEquation_Diag_Uniform}--\eqref{Output_Diag_Uniform} is presented and analyzed. More particularly, the class \eqref{StateEquation_Diag_Uniform}--\eqref{Output_Diag_Uniform} has already been considered in several works, dealing with a wide range of different applications. For instance, Lyapunov exponential stability for \eqref{StateEquation_Diag_Uniform}--\eqref{Output_Diag_Uniform} has been treated in \cite{Bastin_2007, Coron_2007_Lyapunov, Diagne_2012} among others. Therein, the function $\lambda_0$ is in general replaced by a diagonal matrix with possibly different entries. Moreover, the boundary conditions, the inputs and the outputs are slightly different compared to \eqref{Input_Diag_Uniform}--\eqref{Output_Diag_Uniform}. Backstepping for boundary control and output feedback stabilization for networks of linear one-dimensional hyperbolic PDEs has been considered in \cite{Auriol_2020, Auriol_BreschPietri, Vazquez_2011_Backstepping}. Boundary feedback control for systems written as in \eqref{StateEquation_Diag_Uniform}--\eqref{Output_Diag_Uniform} has been tackled in \cite{DeHalleux_2003,Prieur_2008_ConservationLaws,Prieur_Winkin_2018}. In these references, the authors considered quite often systems driven by the Saint-Venant-Exner equations as a particular application for illustrating the results. Questions like controllability and finite-time boundary control have also been studied in \cite{Chitour_2023, Coron_2021_NullContr} and \cite{Auriol_2016,Coron_2021_Stab_FiniteTime} to cite a few. Boundary output feedback stabilization for the class \eqref{StateEquation_Diag_Uniform}--\eqref{Output_Diag_Uniform} with slightly different assumptions has been investigated in \cite{Tanwani2018} in the case of measurement errors, leading to the study of input-to-state stability for those systems. Hyperbolic PDEs similar to \eqref{StateEquation_Diag_Uniform}--\eqref{Output_Diag_Uniform} have attracted attention from a numerical point of view in \cite{Gottlich2017} in which the author proposes a discretization for such systems and then studies the decay rate of the norm of the solutions thanks to this discretization. More recently, the Riesz-spectral property of the homogeneous part (i.e. $u\equiv 0$) of \eqref{StateEquation_Diag_Uniform}--\eqref{Input_Diag_Uniform} has been studied in \cite{HastirJacobZwart_Riesz} where explicit formulas for the eigenvalues and the (generalized) eigenfunctions are depicted. The class \eqref{StateEquation_Diag_Uniform}--\eqref{Output_Diag_Uniform} has also been considered in \cite{HastirJacobZwart_Hinf} where $H^\infty$-controllers are computed, based on the equivalent representation of the system as an infinite-dimensional discrete-time system.
As a general comment, the function $\lambda_0$ is quite often replaced by a diagonal matrix and different versions of the boundary conditions are considered in the aforementioned references. 

When $M\equiv 0$, the class \eqref{StateEquation_Diag_Uniform}--\eqref{Output_Diag_Uniform} fits also the port-Hamiltonian formalism discussed and analyzed in details in \cite{JacobZwart, Villegas} and \cite{Zwart_ESAIM}. In this case, the analysis of the zero dynamics for the class \eqref{StateEquation_Diag_Uniform}--\eqref{Output_Diag_Uniform} has been studied in \cite{JacobMorrisZwart_zeroDynamics}. 

In this work, we describe well-posedness of \eqref{StateEquation_Diag_Uniform}--\eqref{Output_Diag_Uniform} thanks to the class of linear one-dimensional port-Hamiltonian systems analyzed in \cite{JacobZwart}. Next we focus on LQ optimal control. This problem is set up and solved by transforming the original PDEs into a discrete-time infinite-dimensional system with matrices as operator dynamics. The optimal control input is then interpreted in terms of the continuous-time variables, leading to the optimal boundary control. As first example a network of vibrating strings with boundary control forces and boundary measurement velocities is considered. As a second illustration, we consider a co-current heat exchanger heated at the inlet of its external tube, for which the temperature difference between the outlet and the inlet is defined as output.

This manuscript is organized as follows: the well-posedness of the class of considered hyperbolic PDEs is analyzed in Section \ref{sec:WP_Hyperbolic_PDEs}. The considered LQ optimal control problem is described and solved in discrete-time in Section \ref{sec:LQ}. Two different examples are given in Section \ref{sec:examples}, while the conclusions are presented in Section \ref{sec:ccl}.

\section{Well-posedness analysis}
\label{sec:WP_Hyperbolic_PDEs}
An important question for infinite-dimensional boundary controlled, boundary observed systems is the well-posedness. This concept has to be understood as follows: for every $\tilde{z}_0\in X$ and every input $u\in\L^2_{\text{loc}}(0,\infty;\mathbb{C}^p)$, the unique mild solution of \eqref{StateEquation_Diag_Uniform}--\eqref{Input_Diag_Uniform} exists such that the state $\tilde{z}$ and the output $y$ are in the spaces $X$ and $\L^2_{\text{loc}}(0,\infty;\mathbb{C}^m)$, respectively.

Studying well-posedness for \eqref{StateEquation_Diag_Uniform}--\eqref{Output_Diag_Uniform} is eased a lot by considering the system as a boundary controlled port-Hamiltonian system. Before going into the well-posedness part, let us remark that the homogeneous part of \eqref{StateEquation_Diag_Uniform}--\eqref{Input_Diag_Uniform} (i.e. with $u\equiv 0$) may be written as $\dot{\tilde{z}}(t) = A\tilde{z}(t), \tilde{z}(0) = \tilde{z}_0$ with
\begin{align*}
Af &= -\frac{\d}{\d\zeta}(\lambda_0 f) + Mf,\\
D(A) &= \{f\in X \mid \lambda_0 f\in\H^1(0,1;\mathbb{C}^n), \lambda_0(0)Kf(0) + \lambda_0(1)Lf(1) = 0\},
\end{align*}
where the state space $X$ is equipped with the inner product $\langle f,g\rangle_X = \int_0^1 g(\zeta)^*f(\zeta)\d\zeta$. 

To address well-posedness, we start with the following assumption.
\begin{assumption}\label{Assum_M}
The matrix $M$ satisfies $M(\cdot)\in\L^\infty(0,1;\mathbb{R}^{n\times n})$.
\end{assumption}
As a consequence of Assumption \ref{Assum_M}, we have the following lemma.
\begin{lemma}
Under Assumption \ref{Assum_M}, the two following differential equations 
\begin{align}
Q'(\zeta) &= -\frac{1}{\lambda_0(\zeta)}Q(\zeta)M(\zeta),\label{Equa_Diff_Q}\\
P'(\zeta) &= \frac{1}{\lambda_0(\zeta)}M(\zeta)P(\zeta),\label{Equa_Diff_P}
\end{align}
with $\zeta\in (0,1]$ and $Q(0) = I = P(0)$, possess unique continuous and continuously differentiable solutions. Moreover, $P$ and $Q$ are invertible and satisfy 
\begin{equation}
Q(\zeta) = P(\zeta)^{-1}
\label{Link_P_Q}
\end{equation}
for all $\zeta\in [0,1]$.
\end{lemma}
\begin{proof}
Existence and uniqueness of a continuous and continuously differentiable solution of \eqref{Equa_Diff_Q} and \eqref{Equa_Diff_P} follows from \cite{trostorff2023characterisation}. In order to establish \eqref{Link_P_Q}, observe that 
\begin{align*}
\left(Q(\zeta)P(\zeta)\right)' &= Q'(\zeta)P(\zeta) + Q(\zeta)P'(\zeta)\\
&= -\frac{1}{\lambda_0(\zeta)}Q(\zeta)M(\zeta)P(\zeta) + \frac{1}{\lambda_0(\zeta)}Q(\zeta)M(\zeta)P(\zeta)\\
&= 0,
\end{align*}
where \eqref{Equa_Diff_Q} and \eqref{Equa_Diff_P} have been used. Noting that $Q(0)P(0) = I$ concludes the proof.
\end{proof}
In the following lemma, we introduce a change of variables that removes the term $M(\cdot)I$ in the operator $A$.
\begin{lemma}\label{Lemma_Chg_Var}
Let Assumption \ref{Assum_M} hold. Let the multiplication operator $\mathcal{Q}:X\to X$ be defined by $(\mathcal{Q}f)(\zeta) = Q(\zeta)f(\zeta)$ where $Q$ is the unique solution of \eqref{Equa_Diff_Q}. The operator $\mathcal{Q}$ is such that 
\begin{align*}
\mathcal{Q}A\mathcal{Q}^{-1}g &= -\frac{\d}{\d\zeta}(\lambda_0 g),\\
D(\mathcal{Q}A\mathcal{Q}^{-1}) &= \{g\in X \mid \lambda_0 g\in\H^1(0,1;\mathbb{C}^n), \lambda_0(0)Kg(0) + \lambda_0(1)LQ^{-1}(1)g(1) = 0\}.
\end{align*}
\end{lemma}
\begin{proof}
Let us consider $f\in D(A)$ and define $g(\zeta) := Q(\zeta)f(\zeta)$. Hence, $\mathcal{Q}^{-1}g\in D(A)$ and
\begin{align*}
(A\mathcal{Q}^{-1}g)(\zeta) &= -\frac{\d}{\d\zeta}(\lambda_0(\zeta)P(\zeta)g(\zeta)) + M(\zeta)P(\zeta)g(\zeta)\\
&= -\lambda_0'(\zeta)P(\zeta)g(\zeta) - \lambda_0(\zeta)P'(\zeta)g(\zeta) - \lambda_0(\zeta)P(\zeta)g'(\zeta) + M(\zeta)P(\zeta)g(\zeta).
\end{align*}
Using the differential equation \eqref{Equa_Diff_P} this implies that 
\begin{equation*}
(A\mathcal{Q}^{-1}g)(\zeta) = -P(\zeta)\frac{\d}{\d\zeta}(\lambda_0(\zeta)g(\zeta)).
\end{equation*}
As a consequence, there holds $(\mathcal{Q}A\mathcal{Q}^{-1}g)(\zeta) = -\frac{\d}{\d\zeta}(\lambda_0(\zeta)g(\zeta))$. To characterize the boundary conditions satisfied by $g$, observe that $f\in D(A)$ gives $\lambda_0(0)Kf(0) + \lambda_0(1)Lf(1) = 0$. Since $f(\zeta) = P(\zeta)g(\zeta)$, we conclude the proof by noting that
\begin{displaymath}
\lambda_0(0)Kg(0) + \lambda_0(1)LQ^{-1}(1)g(1) = 0.
\end{displaymath}
\end{proof}

\begin{remark}
In the case where the matrix $M$ is constant, it is easy to see that the matrices $Q$ and $P$ solutions of the matrix differential equations \eqref{Equa_Diff_Q} and \eqref{Equa_Diff_P} are given by $Q(\zeta) = e^{-M\int_0^\zeta\lambda_0(\eta)^{-1}\d\eta}$ and $P(\zeta) = e^{M\int_0^\zeta\lambda_0(\eta)^{-1}\d\eta}$, respectively. When $M\equiv 0$, $Q$ and $P$ are both the identity matrix.
\end{remark}

Lemma \ref{Lemma_Chg_Var} enables to describe \eqref{StateEquation_Diag_Uniform}--\eqref{Output_Diag_Uniform} equivalently as the following boundary controlled boundary observed system
\begin{align}
\frac{\partial z}{\partial t}(\zeta,t) &= -\frac{\partial}{\partial \zeta}(\lambda_0(\zeta)z(\zeta,t)),\label{StateEquation_Diag_Uniform_Chg_Var}\\
z(\zeta,0) &= z_0(\zeta),\hspace{0.5cm} \zeta\in (0,1),\label{InitEquation_Diag_Uniform_Chg_Var}\\
\left[\begin{matrix}0\\ I\end{matrix}\right]u(t) &= -\lambda_0(0) Kz(0,t) - \lambda_0(1) L P(1)z(1,t),\label{Input_Diag_Uniform_Chg_Var}\\
y(t) &= -\lambda_0(0) K_y z(0,t) - \lambda_0(1) L_y P(1) z(1,t),\label{Output_Diag_Uniform_Chg_Var}
\end{align}
in which the state $z(\zeta,t)$ is related to $\tilde{z}(\zeta,t)$ in \eqref{StateEquation_Diag_Uniform}--\eqref{Output_Diag_Uniform} by $z(\zeta,t) = Q(\zeta)\tilde{z}(\zeta,t)$. In \eqref{StateEquation_Diag_Uniform_Chg_Var}--\eqref{Output_Diag_Uniform_Chg_Var}, the inputs and the outputs functions $u(t)$ and $y(t)$ remain unchanged compared to \eqref{StateEquation_Diag_Uniform}--\eqref{Output_Diag_Uniform}, they are expressed in terms of $z$ instead of $\tilde{z}$.

As is shown in \cite[Chapter 13]{JacobZwart}, well-posedness of \eqref{StateEquation_Diag_Uniform}--\eqref{Output_Diag_Uniform} (equivalently well-posedness of \eqref{StateEquation_Diag_Uniform_Chg_Var}--\eqref{Output_Diag_Uniform_Chg_Var}) is based on the ability of the operator $A$ (equivalently $\mathcal{Q}A\mathcal{Q}^{-1}$) to be the generator of a $C_0$-semigroup of bounded linear operators. Necessary and sufficient conditions are given in the next proposition.
\begin{prop}\label{Prop_Well_Posed}
The boundary control system \eqref{StateEquation_Diag_Uniform}--\eqref{Output_Diag_Uniform} is well-posed if and only if the matrix $K$ is invertible. 
\end{prop}

\begin{proof}
As is shown in \cite[Chapter 13]{JacobZwart}, well-posedness of \eqref{StateEquation_Diag_Uniform}--\eqref{Output_Diag_Uniform} is equivalent to the fact that the operator $A$, or equivalently the operator $\mathcal{Q}A\mathcal{Q}^{-1}$, is the infinitesimal generator of a $C_0$-semigroup. To characterize this, we rely on \cite[Theorem 1.5]{JacobMorrisZwart}. Therefore, observe that the homogeneous part of \eqref{StateEquation_Diag_Uniform_Chg_Var}--\eqref{Input_Diag_Uniform_Chg_Var} may be written as the following port-Hamiltonian system
\begin{equation}
\left\{\begin{array}{l}
\displaystyle\frac{\partial \mathbf{z}}{\partial t}(\zeta,t) = P_1\frac{\partial}{\partial\zeta}(\mathcal{H}(\zeta)\mathbf{z}(\zeta,t)),\hspace{0.8cm} \zeta\in (0,1),t\geq 0,\\
\mathbf{z}(\zeta,0) = \mathbf{z}_0(\zeta),\hspace{0.8cm} \zeta\in (0,1),\\
0 = \tilde{W}_{B}\left(\begin{matrix}
(\mathcal{H}\mathbf{z})(1,t)\\
(\mathcal{H}\mathbf{z})(0,t)
\end{matrix}\right), \hspace{0.8cm} t\geq 0,
\end{array}\right.
\label{Class_PHS}
\end{equation}
where $P_1 = -I, \mathcal{H}(\zeta) = \lambda_0(\zeta)$ and $\tilde{W}_{B} := \left(\begin{matrix}W_1 & W_0\end{matrix}\right)$ with $W_1 = LP(1)$ and $W_0 = K$. The associated state space is $\mathcal{X} := \L^2(0,1;\mathbb{C}^n)$ equipped with the inner product $\langle f,g\rangle_\mathcal{X} = \int_0^1 g^*(\zeta)\mathcal{H}(\zeta)f(\zeta)\d\zeta$. Following the notations of \cite[Theorem 1.5]{JacobMorrisZwart}, we denote by $Z^+(\zeta)$ the span of eigenvectors of $P_1\mathcal{H}(\zeta)$ corresponding to positive eigenvalues and $Z^-(\zeta)$ the span of eigenvectors corresponding to negative eigenvalues. It is then clear that $Z^+(\zeta) = \{0\}$ and that $Z^-(\zeta) = \mathbb{C}^n$ since the scalar-valued function $\lambda_0$ is positive. Then Theorem 1.5 of \cite{JacobMorrisZwart} states that the operator $\mathcal{Q}A\mathcal{Q}^{-1}$ is the infinitesimal generator of a $C_0$-semigroup if and only if $W_1\mathcal{H}(1)Z^+(1)\oplus W_0\mathcal{H}(0)Z^-(0) = \mathbb{C}^n$. Equivalently, the last equality is $\mathbb{C}^n = L\lambda_0(1)P(1)\{0\}\oplus K\lambda_0(0)\mathbb{C}^n = K\lambda_0(0)\mathbb{C}^n$. Since $\lambda_0(0)$ is not zero, this holds if and only if $K$ is invertible, which concludes the proof.
\end{proof}

\begin{remark}
\begin{itemize} 
\item The class of port-Hamiltonian \eqref{Class_PHS} has been studied in detail in e.g. \cite{JacobMorrisZwart}, \cite{JacobZwart}, \cite{JacobZwart_GAMM}, \cite{Villegas}, where concepts like well-posedness, stability, transfer functions, etc. have been analyzed;
\item The class of systems \eqref{StateEquation_Diag_Uniform_Chg_Var}--\eqref{Output_Diag_Uniform_Chg_Var} may be derived from a port-Hamiltonian system of the form 
\begin{align}
\frac{\partial \mathbf{z}}{\partial t}(\zeta,t) &= \frac{\partial}{\partial \zeta}(\Delta(\zeta)\mathbf{z}(\zeta,t)),\label{StateEquation_Diag}\\
\mathbf{z}(\zeta,0) &= \mathbf{z}_0(\zeta),\hspace{0.8cm} \zeta\in (0,1),\label{InitEquation_Diag}\\
\left[\begin{smallmatrix}0\\ I\end{smallmatrix}\right]u(t) &= \left[\begin{smallmatrix}K_{0+} & K_{0-}\\ K_{u+} & K_{u-}\end{smallmatrix}\right]\left(\begin{smallmatrix}\Lambda(1)\mathbf{z}_+(1,t)\\\Theta(0)\mathbf{z}_-(0,t)\end{smallmatrix}\right) + \left[\begin{smallmatrix}L_{0+} & L_{0-}\\ L_{u+} & L_{u-}\end{smallmatrix}\right]\left(\begin{smallmatrix}\Lambda(0)\mathbf{z}_+(0,t)\\ \Theta(1)\mathbf{z}_-(1,t)\end{smallmatrix}\right),\label{Input_Diag}\\
y(t) &= \left[\begin{smallmatrix}K_{y+} & K_{y-}\end{smallmatrix}\right]\left(\begin{smallmatrix}\Lambda(1)\mathbf{z}_+(1,t)\\\Theta(0)\mathbf{z}_-(0,t)\end{smallmatrix}\right) + \left[\begin{smallmatrix}L_{y+} & L_{y-}\end{smallmatrix}\right]\left(\begin{smallmatrix}\Lambda(0)\mathbf{z}_+(0,t)\\ \Theta(1)\mathbf{z}_-(1,t)\end{smallmatrix}\right)\label{Output_Diag}.
\end{align}
Note that the PDE in \eqref{StateEquation_Diag} corresponds to \eqref{Class_PHS} in which $P_1\mathcal{H}(\zeta)$ is a diagonal matrix which is decomposed as $P_1\mathcal{H}(\zeta) = \left(\begin{smallmatrix}\Lambda(\zeta) & 0\\ 0 & \Theta(\zeta)\end{smallmatrix}\right) =: \Delta(\zeta)$, where $\Lambda(\zeta)\in\mathbb{R}^{k\times k}$ is a diagonal matrix whose diagonal elements are positive functions and $\Theta(\zeta)\in\mathbb{R}^{(n-k)\times (n-k)}$ is a diagonal matrix whose diagonal elements are negative functions. In \eqref{StateEquation_Diag}, the state vector is expressed as $\mathbf{z}(\zeta,t) = \left[\begin{smallmatrix}\mathbf{z}_+(\zeta,t)\\ \mathbf{z}_-(\zeta,t)\end{smallmatrix}\right], \mathbf{z}_+\in\mathbb{C}^k, \mathbf{z}_-\in\mathbb{C}^{n-k}$. In the case where $\Delta$ is constant and the elements on its diagonal are commensurate, it is shown in \cite[Section 4]{JacobMorrisZwart_zeroDynamics} that \eqref{StateEquation_Diag}--\eqref{Output_Diag} may be written as \eqref{StateEquation_Diag_Uniform_Chg_Var}--\eqref{Output_Diag_Uniform_Chg_Var} with $\lambda_0$ constant. In this case the intervals are divided in a series of intervals to adjust the propagation periods. Moreover, the change of variables $\tilde{\mathbf{z}}_l(\zeta,t) := \mathbf{z}_l(1-\zeta,t)$ transforms positive wave speeds into negative ones for the $l$-th component of $\mathbf{z}$. The whole process is detailed in \cite{JacobMorrisZwart_zeroDynamics} and \cite{Suzuki}. This procedure is still valid when $\Delta$ is spatially dependent. In such a situation, the key assumption is that $\eta_j := \int_0^1 \Delta_{jj}^{-1}(\tau)\d\tau$ are commensurate, where $\Delta_{jj}$ is the $j$-th element of the diagonal of $\Delta$.
\end{itemize}
\end{remark}

\section{LQ-optimal control approach}
\label{sec:LQ}
\subsection{Control objective}
For the system \eqref{StateEquation_Diag_Uniform}--\eqref{Output_Diag_Uniform}, our control objective consists in computing the optimal control $u^\text{opt}(\cdot)\in\L^2(0,\infty;\mathbb{C}^p)$ which minimizes the cost 
\begin{equation}
J(z_0,u) = \int_0^\infty \left[\Vert y(t)\Vert^2_{\mathbb{C}^m} + \Vert u(t)\Vert^2_{\mathbb{C}^p}\right]\d t.
\label{Cost_Conti}
\end{equation}
This problem is known as the Linear-Quadratic (LQ) optimal control problem. When the control and observation operators are bounded, it is shown that the solution of the LQ-optimal control problem can be expressed via the unique positive self-adjoint solution of the operator Riccati equation. In this case, the optimal control input is given by the state feedback $\tilde{u}(t) = -B^*Qz(t), t\geq 0$, $Q$ and $B$ being the solution of the aforementioned Riccati equation and the control operator, respectively. The closed-loop system is obtained by connecting the optimal control to the open-loop system. The classical Riccati equation is not valid anymore when the input and the output operators are unbounded, it has then to be adapted. In particular, in \cite{Weiss_Zwart} and \cite{WeissWeiss}, it is mentioned that under the weak regularity of the system and some stability assumptions, the optimal control input is again given by a state feedback. The corresponding feedback operator, and in particular its domain, are often difficult to obtain, except for very particular cases. 

We will see in this section that the computation of the optimal control for \eqref{StateEquation_Diag_Uniform}--\eqref{Output_Diag_Uniform} is simplified a lot, actually it does not require to use the theory of LQ-optimal control with unbounded control and observation operators.

\subsection{Discrete-time systems approach for LQ-optimal control}

Here we link the solution of the LQ-optimal control problem to the solution of an equivalent LQ-optimal control problem in discrete-time. An equivalent representation of \eqref{StateEquation_Diag_Uniform_Chg_Var}--\eqref{Output_Diag_Uniform_Chg_Var} as a discrete-time system is established next.
\begin{lemma}\label{Lemma_Representation}
Let us assume that the matrix $K$ is invertible. Moreover, let the functions $k:[0,1]\to [0,1]$ and $p:[0,1]\to [0,\infty)$ be defined as
\begin{align}
k(\zeta) := 1-p(\zeta)p(1)^{-1}, \hspace{0.8cm} p(\zeta) := \int_0^\zeta \lambda_0(\eta)^{-1}\d\eta,
\label{Functions_k_p}
\end{align}
respectively. In addition, let the matrices $A_d\in\mathbb{C}^{n\times n}, B_d\in\mathbb{C}^{n\times p}, C_d\in\mathbb{C}^{m\times n}$ and $D_d\in\mathbb{C}^{m\times p}$ be respectively given by
\begin{align}
A_d &:= -K^{-1}LP(1),\hspace{0.5cm} B_d := -K^{-1}\left[\begin{matrix}0\\ I\end{matrix}\right],\nonumber\\
C_d &:= (K_y K^{-1}L - L_y)P(1),\hspace{0.5cm} D_d := K_y K^{-1}\left[\begin{matrix}0\\ I\end{matrix}\right].\label{Matrices_Discrete}
\end{align}
The continuous-time system \eqref{StateEquation_Diag_Uniform_Chg_Var}--\eqref{Output_Diag_Uniform_Chg_Var} may be equivalently written as the following discrete-time system
\begin{align}
z_d(j+1)(\zeta) &= A_dz_d(j)(\zeta) + B_d u_d(j)(\zeta),\label{State_Discrete}\\
z_d(0)(k(\zeta)) &= \lambda_0(\zeta)z_0(\zeta),\label{Init_Discrete}\\
y_d(j)(\zeta) &= C_dz_d(j)(\zeta) + D_d u_d(j)(\zeta),\label{Output_Discrete}
\end{align}
whose state, input and output spaces are given by 
\begin{align*}
\mathfrak{X} &:= \L^2(0,1;\mathbb{C}^n),\\
\mathcal{U} &:= \L^2(0,1;\mathbb{C}^p),\\
\mathcal{Y} &:= \L^2(0,1;\mathbb{C}^m),
\end{align*}
respectively. In \eqref{State_Discrete}--\eqref{Output_Discrete} $j\in\mathbb{N}, \zeta\in [0,1], z_d(j)\in\mathfrak{X}, u_d(j)\in \mathcal{U}$ and $y_d(j)\in\mathcal{Y}$. In addition, the functions $z_d, u_d$ and $y_d$ are related to the state, the input and the output of the continuous-time system \eqref{StateEquation_Diag_Uniform_Chg_Var}--\eqref{Output_Diag_Uniform_Chg_Var} via the following relations
\begin{align*}
%z_d(0)(k(\zeta)) &= \lambda_0(\zeta)z_0(\zeta),\\
z_d(j)(\zeta) &= f(j+\zeta),\hspace{0.8cm} j\geq 1,\\
u_d(j)(\zeta) &= u((j+\zeta)p(1)),\hspace{0.8cm} j\in\mathbb{N},\\
y_d(j)(\zeta) &= y((j+\zeta)p(1)),\hspace{0.8cm} j\in\mathbb{N},
\end{align*}
where $z(\zeta,t) = \lambda_0(\zeta)^{-1}f\left(k(\zeta) + p(1)^{-1}t\right)$ and $\lambda_0(\zeta)z_0(\zeta) = f(k(\zeta))$
\end{lemma}

\begin{proof}
First define the functions $k:[0,1]\to[0,1]$ and $p:[0,1]\to\mathbb{R}^+$ by
\begin{align*}
k(\zeta) := 1-p(\zeta)p(1)^{-1}, \hspace{0.8cm} p(\zeta) := \int_0^\zeta \lambda_0(\eta)^{-1}\d\eta.
\end{align*}
Note that the function $p$ is a monotonic function satisfying $p(0) = 0$. In addition, the function $k$ has the properties that $k(0) = 1$ and $k(1) = 0$. Observe that the solution of \eqref{StateEquation_Diag_Uniform_Chg_Var} is given by $z(\zeta,t) = \lambda_0(\zeta)^{-1}f\left(k(\zeta) + p(1)^{-1}t\right)$ for some function $f$. In particular, $f\left(k(\zeta)\right) = \lambda_0(\zeta)z_0(\zeta), \zeta\in [0,1]$. Remark that, according to the definition of $p$, there holds $p(\zeta)p(1)^{-1}\in [0,1]$ for every $\zeta\in [0,1]$. By substituting the expression of $z$ into the initial conditions, the boundary conditions and the output equation \eqref{InitEquation_Diag_Uniform_Chg_Var}--\eqref{Output_Diag_Uniform_Chg_Var}, we get
\begin{align*}
f(k(\zeta)) &= \lambda_0(\zeta)z_0(\zeta),\hspace{0.5cm} \zeta\in [0,1],\\
\left[\begin{matrix}0\\ I\end{matrix}\right]u(t) &= -K f(1+p(1)^{-1}t) - L P(1)f(p(1)^{-1}t),\\
y(t) &= -K_y f(1+p(1)^{-1}t) - L_y P(1) f(p(1)^{-1}t).
\end{align*}
Using the invertibility of $K$, we find
\begin{align}
f(k(\zeta)) &= \lambda_0(\zeta)z_0(\zeta),\hspace{0.5cm} \zeta\in [0,1],\label{Init_f}\\
f(1+p(1)^{-1}t) &= -K^{-1}LP(1)f(p(1)^{-1}t) - K^{-1} \left[\begin{matrix}0\\ I\end{matrix}\right]u(t),\hspace{0.5cm} t\geq 0,\label{State_f}\\
y(t) &= (K_y K^{-1}L-L_y)P(1) f(p(1)^{-1}t) + K_y K^{-1}\left[\begin{matrix}0\\ I\end{matrix}\right]u(t)\label{Output_f}.
\end{align}
By defining the matrices $A_d\in\mathbb{C}^{n\times n}, B_d\in\mathbb{C}^{n\times p}, C_d\in\mathbb{C}^{m\times n}$ and $D_d\in\mathbb{C}^{m\times p}$ by 
\begin{align*}
A_d &:= -K^{-1}LP(1),\\
B_d &:= -K^{-1}\left[\begin{matrix}0\\ I\end{matrix}\right],\\
C_d &:= (K_yK^{-1}L-L_y)P(1),\\
D_d &:= K_yK^{-1}\left[\begin{matrix}0\\ I\end{matrix}\right],
\end{align*}
\eqref{State_f}--\eqref{Output_f} may be written equivalently as
\begin{align}
f(1+p(1)^{-1}t) &= A_d f(p(1)^{-1}t) + B_d u(t),\hspace{0.5cm} t\geq 0,\label{State_Equiv_f}\\
y(t) &= C_d f(p(1)^{-1}t) + D_d u(t).\label{Output_Equiv_f}
\end{align}
Now take $j\in\mathbb{N}, \zeta\in [0,1]$ and define
\begin{align*}
z_d(j)\in\L^2(0,1;\mathbb{C}^n),\hspace{0.5cm} u_d(j)\in\L^2(0,1;\mathbb{C}^p),\hspace{0.5cm} y_d(j)\in\L^2(0,1;\mathbb{C}^m)
\end{align*}
by 
\begin{align}
z_d(0)(k(\zeta)) &:= \lambda_0(\zeta)z_0(\zeta),\label{Init_Var}\\
z_d(j)(\zeta) &:= f(j+\zeta),\hspace{0.5cm} j\geq 1,\label{State_Var}\\
u_d(j)(\zeta) &:= u((j+\zeta)p(1)),\hspace{0.5cm} j\in\mathbb{N},\label{Input_Var}\\
y_d(j)(\zeta) &:= y((j+\zeta)p(1)),\hspace{0.5cm} j\in\mathbb{N}.\label{IO_Var}
\end{align}
Remark that, for any $t\geq 0$, we can find $j\in\mathbb{N}$ and $\zeta\in [0,1]$ such that $j+\zeta = p(1)^{-1}t$. With this definition, we get 
\begin{align*}
f(p(1)^{-1}t) &= f(j+\zeta) = z_d(j)(\zeta),\\
u(t) &= u((j+\zeta)p(1)) = u_d(j)(\zeta),\\
y(t) &= y_d(j)(\zeta),
\end{align*}
which allows us to write \eqref{State_Equiv_f} and \eqref{Output_Equiv_f} as \eqref{State_Discrete} and \eqref{Output_Discrete}, respectively.
\end{proof}

\begin{remark}\label{Remark_System_Th_Properties}
\begin{itemize}
\item A procedure similar to the one developed in the proof of Lemma \ref{Lemma_Representation} may be found in \cite[Section 4]{JacobMorrisZwart_zeroDynamics} in the case of constant function $\lambda_0$ and $M\equiv 0$.
\item The idea that is followed in Lemma \ref{Lemma_Representation} can be compared to the one followed in \cite{Yamamoto_Discrete}. Therein, the author looks at finite-dimensional continuous-time systems by cutting the time interval into pieces of same length, giving rise to an equivalent infinite-dimensional discrete-time representation of the original systems.
\item As the function $p$ is monotonic, it is injective and the range of $p(\cdot)p(1)^{-1}$ equals $[0,1]$. This has the consequence that, knowing the initial condition $z_0$ on the interval $[0,1]$ is equivalent in the knowledge of $z_d(0)$ on $[0,1]$.
\item The discrete-time representation \eqref{State_Discrete}--\eqref{Output_Discrete} is interesting for system theoretic properties of the continuous-time system \eqref{StateEquation_Diag_Uniform_Chg_Var}--\eqref{Output_Diag_Uniform_Chg_Var}. For instance, exponential stability of \eqref{StateEquation_Diag_Uniform_Chg_Var}--\eqref{Output_Diag_Uniform_Chg_Var} is equivalent to $r(A_d)<1$, with $r$ being the spectral radius. This condition can also be found in \cite[Theorem 11]{JacobMorrisZwart_zeroDynamics}. 
\end{itemize}
\end{remark}

The following proposition gives an equivalent representation of the cost \eqref{Cost_Conti} with the variables of the discrete-time system.

\begin{prop}\label{Cost_Equiv_Conti_Discrete}
The cost \eqref{Cost_Conti} associated to the optimal control problem for the continuous-time system may be rewritten with the discrete-time input and output variables as
\begin{equation}
J(z_0,u) = p(1)\sum_{j=0}^\infty \left[\Vert y_d(j)(\cdot)\Vert_\mathcal{Y}^2 + \Vert u_d(j)(\cdot)\Vert_\mathcal{U}^2\right] =: p(1)J_d(z_d(0),u_d).
\label{Cost_Discrete}
\end{equation}
\end{prop}

\begin{proof}
Let us take $z_0\in X$ and for $t\geq 0$ denote by $u(t)\in\mathbb{C}^p$ and $y(t)\in\mathbb{C}^m$ the input and the output of \eqref{StateEquation_Diag_Uniform}--\eqref{Output_Diag_Uniform} at time $t$. First observe that
\begin{align*}
J(z_0,u) &= \int_0^\infty \left[\Vert y(t)\Vert_{\mathbb{C}^m}^2 + \Vert u(t)\Vert_{\mathbb{C}^p}^2\right]\d t\\
&= p(1)\int_0^\infty \left[\Vert y(p(1)t)\Vert^2_{\mathbb{C}^m} + \Vert u(p(1)t)\Vert_{\mathbb{C}^m}^2\right]\d t.
\end{align*}
%As the number $p(1)$ is positive, minimizing $J(z_0,u)$ or $p(1)^{-1}J(z_0,u)$ is equivalent. For this reason, we shall keep the cost $p(1)^{-1}J(z_0,u)$ in what follows. 
Moreover, there holds
\begin{align*}
p(1)^{-1}J(z_0,u) &= \int_0^\infty \left[\Vert y(p(1) t)\Vert_{\mathbb{C}^m}^2 + \Vert u(p(1) t)\Vert_{\mathbb{C}^p}^2\right]\d t\\
&= \sum_{j=0}^\infty \int_j^{j+1} \left[\Vert y(p(1)t)\Vert_{\mathbb{C}^m}^2 + \Vert u(p(1)t)\Vert_{\mathbb{C}^p}^2\right]\d t\\
&= \sum_{j=0}^\infty \int_0^1\Vert y(p(1)(\tau+j))\Vert_{\mathbb{C}^m}^2\d \tau + \int_0^1\Vert u(p(1)(\tau+j))\Vert_{\mathbb{C}^p}^2\d\tau.
\end{align*}
Thanks to the relations between the input and the output functions of the continuous-time system \eqref{StateEquation_Diag_Uniform_Chg_Var}--\eqref{Output_Diag_Uniform_Chg_Var} and the discrete-time system \eqref{State_Discrete}--\eqref{Output_Discrete}, namely $u(p(1)(\tau+j)) = u_d(j)(\tau)$ and $y(p(1)(\tau+j)) = y_d(j)(\tau)$, there holds 
\begin{align*}
p(1)^{-1}J(z_0,u) &= \sum_{j=0}^\infty \left[\Vert y_d(j)\Vert^2_\mathcal{Y} + \Vert u_d(j)\Vert^2_\mathcal{U}\right] =: J_d(z_d(0),u_d),
\end{align*}
where the input and the output spaces associated to the discrete-time system are given by $\mathcal{U} = \L^2(0,1;\mathbb{C}^p)$ and $\mathcal{Y} = \L^2(0,1;\mathbb{C}^m)$. This concludes the proof.
\end{proof}

Note that the invertibility of the matrix $K$ is assumed from now on, so that the continuous-time system \eqref{StateEquation_Diag_Uniform_Chg_Var}--\eqref{Output_Diag_Uniform_Chg_Var} is well-posed and so that the matrices $A_d, B_d, C_d$ and $D_d$ are well defined.

The operators in the infinite-dimensional discrete-time system \eqref{State_Discrete}--\eqref{Output_Discrete} are constant multiplication operators represented by matrices, which makes them bounded. Now we take the boundary conditions, the input and the output of the continuous-time system \eqref{Input_Diag_Uniform_Chg_Var}--\eqref{Output_Diag_Uniform_Chg_Var}. In order to emphasize better the matrices $A_d, B_d, C_d$ and $D_d$ in the continuous-time system \eqref{StateEquation_Diag_Uniform_Chg_Var}--\eqref{Output_Diag_Uniform_Chg_Var}, observe that
\begin{align*}
\lambda_0(0)z(0,t) &= -K^{-1}\left[\begin{smallmatrix}0\\ I\end{smallmatrix}\right]u(t) - \lambda_0(1)K^{-1}LP(1)z(1,t) =B_d u(t) + \lambda_0(1)A_dz(1,t).
\end{align*}
Incorporating the previous expression in \eqref{Output_Diag_Uniform_Chg_Var} this implies that
\begin{align*}
y(t) &= -K_yB_du(t) + (-K_yA_d - L_yP(1))\lambda_0(1)z(1,t) = D_d u(t) + C_d\lambda_0(1)z(1,t).
\end{align*}
This means that the continuous-time system \eqref{StateEquation_Diag_Uniform_Chg_Var}--\eqref{Output_Diag_Uniform_Chg_Var} may be written as
\begin{align}
\frac{\partial z}{\partial t}(\zeta,t) &= -\frac{\partial}{\partial \zeta}(\lambda_0(\zeta)z(\zeta,t)),\label{StateEquation_Equiv_Conti}\\
z(\zeta,0) &= z_0(\zeta),\label{InitEquation_Equiv_Conti}\\
B_d u(t) &= \lambda_0(0)z(0,t) - A_d\lambda_0(1)z(1,t),\label{Input_Equiv_Conti}\\
y(t) &= C_d\lambda_0(1)z(1,t) + D_d u(t).\label{Output_Equiv_Conti}
\end{align}
\begin{remark}
The assumptions needed for the resolution of the LQ-optimal control problem for the continuous-time system \eqref{StateEquation_Diag_Uniform_Chg_Var}--\eqref{Output_Diag_Uniform_Chg_Var} (i.e. \eqref{StateEquation_Equiv_Conti}--\eqref{Output_Equiv_Conti}) will be formulated in terms of the matrices $A_d, B_d, C_d$ and $D_d$.
\end{remark}

For \eqref{StateEquation_Equiv_Conti}--\eqref{Output_Equiv_Conti} with the cost $p(1)^{-1}J(z_0,u)$, see \eqref{Cost_Discrete}, we characterize the LQ optimal control in terms of existence, uniqueness and ability to stabilize the open-loop system \eqref{State_Discrete}, relying on \cite[Lemmas 3.3, 3.6 and 3.7]{Opmeer_LQ_Discrete}. 
To this end, let us introduce the two following notions for infinite-dimensional discrete-time systems, see \cite[Section 3]{Opmeer_LQ_Discrete}.
\begin{defn}
Consider $(A_d,B_d,C_d,D_d)$ an infinite-dimensional discrete-time system. By $\mathcal{C}$, we denote the observability map of the discrete-time system, which is defined as $(\mathcal{C}x)_j = C_dA_d^jx, j\in\mathbb{N}$. The discrete-time system is said to be output stable if the image of $\mathcal{C}$ is contained in $l^2(\mathbb{N};\mathcal{Y})$.
\end{defn}

\begin{defn}
Consider $(A_d,B_d,C_d,D_d)$ an infinite-dimensional discrete-time system. The controllability map $\mathcal{B}$ of the discrete-time system is defined for finitely nonzero $\mathcal{U}$-valued sequences $u$ by $\mathcal{B}u := \sum_{j=0}^\infty A^j_d B_d u_{-j-1}$. The discrete-time system is said to be input stable if $\mathcal{B}$ extends to a bounded map from $l^2(\mathbb{Z}^-;\mathcal{U})$ to $\mathfrak{X}$.
\end{defn}

An inherent concept to output and input stability is the notion of \textit{strong stability}. The discrete-time system $(A_d,B_d,C_d,D_d)$ is said to be strongly stable if for all $x\in\mathfrak{X}$, there holds $A_d^n x\to 0$ whenever $n\to\infty$.

\begin{remark}\label{Remark_IS_OS}
Characterizations of output and input stability are given in \cite{Opmeer_LQ_Discrete} in terms of Lyapunov equations. In particular, the discrete-time system is output stable if and only if the observation Lyapunov equation 
\begin{equation}
A_d^*L_d^oA_d - L_d^o + C_d^*C_d = 0
\label{LyapunovEquation_OS}
\end{equation}
has a nonnegative self-adjoint solution $L_d^o\in\mathbb{C}^{n\times n}$, see \cite[Lemma 3.1]{Opmeer_LQ_Discrete}. Input stability is equivalent to the existence of a nonnegative self-adjoint solution $L_d^c\in\mathbb{C}^{n\times n}$ of the following control Lypaunov equation 
\begin{equation}
A_dL_d^cA_d^* - L_d^c + B_dB_d^* = 0,
\label{LyapunovEquation_IS}
\end{equation} 
see \cite[Lemma 3.4]{Opmeer_LQ_Discrete}. As the operators $A_d, B_d, C_d$ and $D_d$ are matrices, input and output stability conditions are finite-dimensional conditions. In particular, equations \eqref{LyapunovEquation_OS} and \eqref{LyapunovEquation_IS} are
\begin{equation*}
P(1)^*L^*K^{-*}L_d^oK^{-1}LP(1) - L_d^o + P(1)^*(L^*K^{-*}K_y^*-L_y^*)(K_yK^{-1}L-L_y)P(1) = 0
\end{equation*}
and 
\begin{equation*}
LP(1) L_d^c P(1)^*L^* - K L_d^c K^* + \left(\begin{smallmatrix}0_{(n-p)\times(n-p)} & 0_{(n-p)\times p}\\ 0_{p\times (n-p)} & I_{p\times p}\end{smallmatrix}\right) = 0.
\end{equation*}
\end{remark}

The concepts of output and input stability are often too strong and are actually not needed for the LQ-optimal control problem to be well-posed. For this reason, let us introduce the concepts of output and input stabilizability, see \cite[Section 3]{Opmeer_LQ_Discrete}.

By output stabilizability, we mean that there exists an operator $F_d\in\mathcal{L}(\mathfrak{X},\mathcal{U})$ such that $(A_d+B_dF_d,0,\left[\begin{smallmatrix} F_d\\ C_d+D_dF_d\end{smallmatrix}\right],0)$ is output stable. The discrete-time system $(A_d,B_d,C_d,D_d)$ is input stabilizable if there exists an operator $L_d\in\mathcal{L}(\mathcal{Y},\mathfrak{X})$ such that $(A_d+L_dC_d,[\begin{smallmatrix}L_d & B_d+L_dD_d\end{smallmatrix}],0,0)$ is input stable.

A key assumption in the resolution of the LQ-optimal control problem is the finite-cost condition, also called optimizability, which means that for every initial condition $z_0\in X$, there exists an input $\tilde{u}\in\L^2(0,\infty;\mathbb{C}^p)$ such that $J(z_0,\tilde{u})<\infty$. For instance, this question has been studied in \cite{JacobZwart_SIAM_Opti} for admissible control operators and bounded output operators. For the continuous-time system \eqref{StateEquation_Equiv_Conti}--\eqref{Output_Equiv_Conti}, we have the following lemma.
\begin{lemma}\label{Optimizability}
The continuous-time system \eqref{StateEquation_Equiv_Conti}--\eqref{Output_Equiv_Conti} with the cost \eqref{Cost_Conti} is optimizable if and only if there exists a nonnegative self-adjoint solution $\Pi\in\mathbb{C}^{n\times n}$ to the following control algebraic Riccati equation (CARE)
\begin{align}
&A_d^*\Pi A_d - \Pi + C_d^*C_d\nonumber\\
&\hspace{2cm}= (C_d^*D_d+A_d^*\Pi B_d)(I+D_d^*D_d + B_d^*\Pi B_d)^{-1}(D_d^*C_d + B_d^*\Pi A_d).
\label{CARE}
\end{align}
\end{lemma}
\begin{proof}
First observe that optimizability of \eqref{StateEquation_Equiv_Conti}--\eqref{Output_Equiv_Conti} with the cost \eqref{Cost_Conti} is equivalent to the finite-cost condition. Then, by Proposition \ref{Cost_Equiv_Conti_Discrete}, we get that the costs of the continuous-time and the discrete-time systems are related by $J(z_0,u) = p(1)J_d(z_d(0),u_d)$. Thanks to Lemma \ref{Lemma_Representation}, the discrete-time system \eqref{State_Discrete}--\eqref{Output_Discrete} is an equivalent representation of the continuous-time system \eqref{StateEquation_Equiv_Conti}--\eqref{Output_Equiv_Conti}. Hence, $J(z_0,u)<\infty$ if and only if $J_d(z_d(0),u_d)<\infty$. Equivalently, the discrete-time system \eqref{State_Discrete}--\eqref{Output_Discrete} with the cost $J_d$ has to satisfy the finite-cost condition. By \cite[Lemma 3.3]{Opmeer_LQ_Discrete}, this is equivalent to the existence of a nonnegative self-adjoint solution to the CARE \eqref{CARE}.
\end{proof}

Sufficient conditions that guarantee the uniqueness of the solution $\Pi\in\mathbb{C}^{n\times n}$ of the CARE \eqref{CARE} are given in the following lemma.
\begin{lemma}
Let us suppose that the CARE \eqref{CARE} has a nonnegative self-adjoint solution $\Pi\in\mathbb{C}^{n\times n}$. In addition, assume that the Filter Algebraic Riccati Equation (FARE)
\begin{align}
&A_d\tilde{\Pi} A_d^* - \tilde{\Pi} + B_d^*B_d\nonumber\\
&\hspace{2cm}= (B_dD_d^* + A_d\tilde{\Pi} C_d^*)(I + D_dD_d^* + C_d\tilde{\Pi} C_d^*)^{-1}(D_d B_d^* + C_d\tilde{\Pi} A_d^*).
\label{FARE}
\end{align}
has a nonnegative self-adjoint solution $\tilde{\Pi}\in\mathbb{C}^{n\times n}$. Moreover, assume that the matrix 
\begin{equation}
A_\Pi := A_d - B_d(I + D_d^*D_d + B_d^*\Pi B_d)^{-1}(D_d^*C_d + B^*_d\Pi A_d)
\label{A_Pi}
\end{equation}
is stable in the sense that $r(A_\Pi)<1$. Then $\Pi$ is the unique nonnegative solution of \eqref{CARE}.
\end{lemma}
\begin{proof}
See \cite[Lemma 3.7]{Opmeer_LQ_Discrete}.
\end{proof}
\begin{remark}\label{Charact_IO_Stabilizab}
Existence of nonnegative self-adjoint solutions to the CARE \eqref{CARE} and to the FARE \eqref{FARE} are also called output and input stabilizability of the discrete-time system $(A_d, B_d, C_d, D_d)$, see \cite[Lemmas 3.3 and 3.6]{Opmeer_LQ_Discrete}. One main advantage of dealing with the discrete-time representation of \eqref{StateEquation_Equiv_Conti}--\eqref{Output_Equiv_Conti} is that the operators $A_d, B_d, C_d$ and $D_d$ are matrices, which entails that the equations \eqref{CARE} and \eqref{FARE} are finite-dimensional equations.
\end{remark}
Observe now that, thanks to Lemma \ref{Optimizability}, the existence of a nonnegative self-adjoint solution of the FARE \eqref{FARE} is equivalent to the optimizability of the dual system $(A_d^*, C_d^*, B_d^*, D_d^*)$ in discrete-time. 

The following theorem characterizes the solution of the LQ-optimal control problem for the continuous-time system \eqref{StateEquation_Equiv_Conti}--\eqref{Output_Equiv_Conti} with the cost $p(1)^{-1}J(z_0,u)$, see \eqref{Cost_Discrete}.
\begin{thm}\label{Thm_Unique_Sol}
Let us suppose that the CARE \eqref{CARE} and the FARE \eqref{FARE} have nonnegative self-adjoint solutions $\Pi\in\mathbb{C}^{n\times n}$ and $\tilde{\Pi}\in\mathbb{C}^{n\times n}$, respectively. Moreover, assume that the matrix $A_\Pi$ given in \eqref{A_Pi} satisfies $r(A_\Pi)<1$. Then the unique optimal control $u^{\text{opt}}(t)$ that minimizes the cost functional $p(1)^{-1}J(z_0,u)$, see \eqref{Cost_Conti}, is given by 
\begin{equation}
u^{\text{opt}}(t) = \lambda_0(1)F_dz^{\text{opt}}(1,t),
\label{Optimal_Input_Original_Var}
\end{equation}
where $z^{\text{opt}}(\cdot,t)$ is the state trajectory of the closed-loop system 
\begin{align}
\frac{\partial z^{\text{opt}}}{\partial t}(\zeta,t) &= -\frac{\partial}{\partial \zeta}(\lambda_0(\zeta)z^\text{opt}(\zeta,t)),\label{State_Closed-loop}\\
z^\text{opt}(\zeta,0) &= z_0(\zeta),\label{Init_Closed-loop}\\
\lambda_0(0)z^{\text{opt}}(0,t) &= (A_d + B_dF_d)\lambda_0(1)z^\text{opt}(1,t),\label{BC_Closed-loop}
\end{align}
in which $F_d := -(I + D_d^*D_d + B_d^*\Pi B_d)^{-1}(D_d^*C_d + B^*_d\Pi A_d)$. Moreover, the closed-loop system \eqref{State_Closed-loop}--\eqref{BC_Closed-loop} is exponentially stable.
\end{thm}
\begin{proof}
To establish \eqref{Optimal_Input_Original_Var}, first note that the optimal control in discrete-time, denoted by $u_d^{\text{opt}}(j)$, is given by $u_d^{\text{opt}}(j) = F_dz_d^{\text{opt}}(j)$ where $z_d^\text{opt}$ is the state-trajectory of the closed-loop system in discrete-time\footnote{By closed-loop system in discrete-time, we mean the discrete-time system \eqref{State_Discrete}--\eqref{Output_Discrete} in closed-loop with the optimal control input $u_d^{\text{opt}}(j) = F_dz_d^{\text{opt}}(j)$.}, that is,
\begin{align*}
z_d^\text{opt}(j+1)(\zeta) &= A_\Pi z_d^\text{opt}(j)(\zeta),\\
z_d^\text{opt}(0)(k(\zeta)) &= \lambda_0(\zeta)z_0(\zeta).
\end{align*}
As is performed in Lemma \ref{Lemma_Representation}, the trajectory $z_d^\text{opt}(j)(\zeta)$ corresponds to a function $f^\text{opt}(j+\zeta)$ that satisfies $\lambda_0(\zeta)^{-1} f^{\text{opt}}(k(\zeta)+p(1)^{-1}t) = z^{\text{opt}}(\zeta,t)$, where $z^\text{opt}(\cdot,t)$ is the optimal state trajectory in continuous-time. Then take $j\in\mathbb{N}, \zeta\in [0,1]$ and $t\geq 0$ such that $j+\zeta = p(1)^{-1}t$. There holds
\begin{equation}
u^{\text{opt}}(t) = u_d^{\text{opt}}(j)(\zeta) = F_d f^{\text{opt}}(j+\zeta) = F_d f^{\text{opt}}(p(1)^{-1}t) = \lambda_0(1)F_d z^\text{opt}(1,t).
\end{equation}
Attaching the input $u^\text{opt}$ to the system \eqref{StateEquation_Equiv_Conti}--\eqref{Input_Equiv_Conti} implies that the closed-loop system in continuous-time is given by \eqref{State_Closed-loop}--\eqref{BC_Closed-loop}. According to Remark \ref{Remark_System_Th_Properties}, it is exponentially stable if and only if $r(A_d + B_dF_d)<1$, which concludes the proof.
\end{proof}

\begin{remark}\label{Rem:Scaling}
\begin{itemize}
\item If the cost associated to the LQ-optimal control problem is given by $p(1)J_d$, then the solution of the associated CARE is given by $\hat{\Pi} = p(1)\Pi$. Indeed, by using the notation $p(1) =: \lambda$, the cost $p(1)J_d$ may be written as
\begin{align*}
\lambda J_d = \sum_{j=0}^\infty \left[\Vert \tilde{C}_d z_d(j)(\cdot) + D_d\tilde{u}_d(j)(\cdot)\Vert_\mathcal{Y}^2 + \Vert \tilde{u}_d(j)(\cdot)\Vert^2_\mathcal{U}\right],
\end{align*}
where $\tilde{C}_d := \lambda^{\frac{1}{2}}C_d$ and $\tilde{u}_d(j) := \lambda^\frac{1}{2}u_d(j)$. The state equation \eqref{State_Discrete} is then equivalent to $z_d(j+1) = A_dz_d(j) + \tilde{B}_d\tilde{u}_d(j)$ with $\tilde{B}_d := \lambda^{-\frac{1}{2}}B_d$. It is then easy to see that $\hat{\Pi} = \lambda\Pi$ solves the CARE \eqref{CARE} in which $B_d$ and $C_d$ are replaced by $\tilde{B}_d$ and $\tilde{C}_d$ if and only if $\Pi$ solves \eqref{CARE}. The corresponding optimal feedback operator $\tilde{F}_d$ is then given by $\tilde{F}_d = \lambda^{\frac{1}{2}}F_d$ while the closed-loop operator satisfies $A_d + \tilde{B}_d\tilde{F}_d = A_d + B_dF_d$.
\item The optimal cost associated to the optimal control is given by $J_d(z_d(0),u_d^{\text{opt}}) = \langle z_d^\text{opt}(0),\Pi z_d^\text{opt}(0)\rangle_\mathcal{X}$.
\item The optimal control is a boundary control which only affects the boundary conditions of \eqref{StateEquation_Equiv_Conti}--\eqref{Output_Equiv_Conti}. The PDE itself remains unchanged.
\end{itemize}
\end{remark}

\section{Examples}
\label{sec:examples}

\subsection{A network of interconnected vibrating strings}
\begin{figure}
\begin{center}
\includegraphics[scale=1]{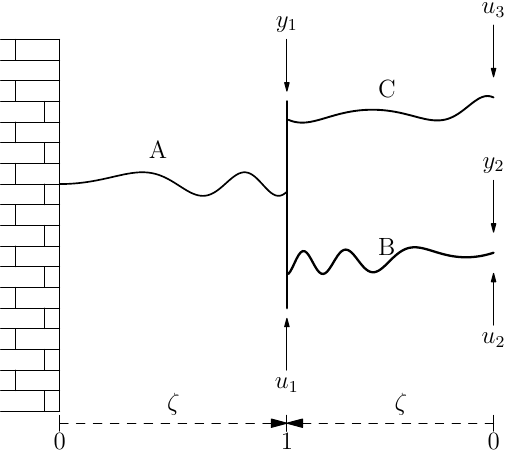}
\caption{Schematic view of three interconnected vibrating strings.\label{fig:ConnectedStrings}}
\end{center}
\end{figure}

We consider a network of three interconnected vibrating strings of same lengths, namely $1$, which are connected with a mass-less bar in the middle, as depicted in Figure \ref{fig:ConnectedStrings}. 
For string $\text{A}$, we let the spatial variable go from $0$ to $1$ whereas it goes in the opposite direction for strings $\text{B}$ and $\text{C}$. We denote by $\rho_k$ and $T_k$ the mass-density and the Young's modulus of the string $k$, respectively, $k\in\{\text{A,B,C}\}$. Each string is described by the following PDE
\begin{equation}
\frac{\partial}{\partial t}\left[\begin{matrix}\rho_k\frac{\partial w_k}{\partial t}\\ \frac{\partial w_k}{\partial\zeta}\end{matrix}\right](\zeta,t) = \left[\begin{matrix}0 & 1\\ 1 & 0\end{matrix}\right]\frac{\partial}{\partial\zeta}\left(\left[\begin{matrix}\frac{1}{\rho_k(\zeta)} & 0\\ 0 & T_k(\zeta)\end{matrix}\right]\left[\begin{matrix}\rho_k\frac{\partial w_k}{\partial t}\\ \frac{\partial w_k}{\partial\zeta}\end{matrix}\right](\zeta,t)\right),
\label{Model_Connected_Strings}
\end{equation}
where $k\in\{\text{A,B,C}\}$. The variable $w_k(\zeta,t)$ represents the vertical displacement of string $k$ at position $\zeta\in[0,1]$ and time $t\geq 0$. Forces are applied at $\zeta=0$ of string $\text{B}$ and $\text{C}$, and on the bar in the middle. Moreover, we assume that the velocity of the bar and that of string $\text{B}$ at $\zeta=0$ are measured. Furthermore, string $\text{A}$ is fixed at $\zeta=0$. These facts lead to the following boundary conditions
\begin{align}
\frac{\partial w_\text{A}}{\partial t}(0,t) &= 0,\hspace{0.3cm} \frac{\partial w_\text{A}}{\partial t}(1,t) = \frac{\partial w_{\text{B}}}{\partial t}(1,t),\hspace{0.3cm} \frac{\partial w_{\text{B}}}{\partial t}(1,t) = \frac{\partial w_{\text{C}}}{\partial t}(1,t),\label{BC}\\
y_1(t) &= \frac{\partial w_{\text{C}}}{\partial t}(1,t),\hspace{0.3cm} y_2(t) = \frac{\partial w_{\text{B}}}{\partial t}(0,t),\label{Outputs}\\
u_1(t) &= T_\text{A}(1)\frac{\partial w_\text{A}}{\partial\zeta}(1,t) - T_{\text{B}}(1)\frac{\partial w_{\text{B}}}{\partial \zeta}(1,t) - T_{\text{C}}(1)\frac{\partial w_{\text{C}}}{\partial\zeta}(1,t),\label{InputOne}\\
u_2(t) &= T_{\text{B}}(0)\frac{\partial w_{\text{B}}}{\partial\zeta}(0,t),\hspace{0.3cm} u_3(t) = T_{\text{C}}(0)\frac{\partial w_{\text{C}}}{\partial\zeta}(0,t),\label{InputTwo}
\end{align}
where $u_i(t), i=1,2,3,$ and $y_j(t), j=1, 2,$ are the inputs and the outputs, respectively. We define the state variable $x := \left(\begin{smallmatrix}\rho_\text{A}\frac{\partial w_{\text{A}}}{\partial t} & \frac{\partial w_\text{A}}{\partial\zeta} &\rho_\text{B}\frac{\partial w_{\text{B}}}{\partial t} & \frac{\partial w_\text{B}}{\partial\zeta} & \rho_\text{C}\frac{\partial w_{\text{C}}}{\partial t} & \frac{\partial w_\text{C}}{\partial\zeta}\end{smallmatrix}\right)^T$, the input vector $u(t) := \left(\begin{smallmatrix}u_1(t) & u_2(t) & u_3(t)\end{smallmatrix}\right)^T$ and the output vector $y(t) := \left(\begin{smallmatrix}y_1(t) & y_2(t)\end{smallmatrix}\right)^T$. Hence, \eqref{Model_Connected_Strings} may be written as $\frac{\partial x}{\partial t}(\zeta,t) = P_1\frac{\partial}{\partial \zeta}(\mathcal{H}x)(\zeta,t)$, with 
\begin{align*}
P_1 &= \text{diag}(P,P,P),\hspace{0.5cm} \mathcal{H} = \text{diag}(\mathcal{H}_\text{A},\mathcal{H}_{\text{B}},\mathcal{H}_{\text{C}}),\\
P &= \left(\begin{smallmatrix}0 & 1\\1 & 0\end{smallmatrix}\right),\hspace{0.5cm} \mathcal{H}_k = \left(\begin{smallmatrix}\frac{1}{\rho_k(\zeta)} & 0\\0 & T_k(\zeta)\end{smallmatrix}\right), k\in\{\text{A, B, C}\}.
\end{align*}
Moreover, \eqref{BC}--\eqref{InputTwo} may be rewritten in the following condensed form
\begin{align*}
0 &= \tilde{W}_{B,1}\left(\begin{matrix}\mathcal{H}x(1,t)\\\mathcal{H}x(0,t)\end{matrix}\right),\hspace{0.5cm} u(t) = \tilde{W}_{B,2}\left(\begin{matrix}\mathcal{H}x(1,t)\\\mathcal{H}x(0,t)\end{matrix}\right),\hspace{0.5cm} y(t) = \tilde{W}_C\left(\begin{matrix}\mathcal{H}x(1,t)\\\mathcal{H}x(0,t)\end{matrix}\right),
\end{align*}
where the matrices $\tilde{W}_{B,1}, \tilde{W}_{B,2}$ and $\tilde{W}_C$ are respectively given by
\begin{align*}
\tilde{W}_{B,1} &= \left(\begin{smallmatrix}0 & 0 & 0 & 0 & 0 & 0 & 1 & 0 & 0 & 0 & 0 & 0\\
1 & 0 & -1 & 0 & 0 & 0 & 0 & 0 & 0 & 0 & 0 & 0\\
0 & 0 & 1 & 0 & -1 & 0 & 0 & 0 & 0 & 0 & 0 & 0
\end{smallmatrix}\right),\\
\tilde{W}_{B,2} &= \left(\begin{smallmatrix}0 & 1 & 0 & -1 & 0 & -1 & 0 & 0 & 0 & 0 & 0 & 0\\
0 & 0 & 0 & 0 & 0 & 0 & 0 & 0 & 0 & 1 & 0 & 0\\
0 & 0 & 0 & 0 & 0 & 0 & 0 & 0 & 0 & 0 & 0 & 1
\end{smallmatrix}\right),\\
\tilde{W}_C &= \left(\begin{smallmatrix}0 & 0 & 0 & 0 & 1 & 0 & 0 & 0 & 0 & 0 & 0 & 0\\
0 & 0 & 0 & 0 & 0 & 0 & 0 & 0 & 1 & 0 & 0 & 0
\end{smallmatrix}\right).
\end{align*}
We assume that the three strings are the same in the sense that $\rho_{\text{A}} = \rho_{\text{B}} = \rho_{\text{C}} =: \rho$ and $T_{\text{A}} = T_{\text{B}} = T_{\text{C}} =: T$. Moreover, we denote by $v$ and $\sigma$ the constants $\sqrt{\frac{T}{\rho}}$ and $\frac{1}{\sqrt{\rho T}}$, respectively. As a consequence, $P_1\mathcal{H}$ may be written as $P_1\mathcal{H} = S^{-1}\Delta S$, where $S = \text{diag}(s,s,s), \Delta = \text{diag}(\delta,\delta,\delta), s = \left(\begin{smallmatrix}\frac{1}{2} & \frac{1}{2\sigma}\\-\frac{1}{2} & \frac{1}{2\sigma}\end{smallmatrix}\right), \delta = \left(\begin{smallmatrix}v & 0\\0 & -v\end{smallmatrix}\right)$. By defining $z(\zeta,t) = Sx(\zeta,t)$, we get that $\frac{\partial z}{\partial t}(\zeta,t) = \Delta\frac{\partial z}{\partial\zeta}(\zeta,t)$ with 
\begin{align*}
\left(\begin{matrix}
0_{3\times 3}\\
I_{3\times 3}
\end{matrix}\right)u(t) &= \left(\begin{smallmatrix}0 & 0 & 0 & 0 & 0 & 0 & \sigma v & -\sigma v & 0 & 0 & 0 & 0\\
\sigma v & -\sigma v & -\sigma v & \sigma v & 0 & 0 & 0 & 0 & 0 & 0 & 0 & 0\\
0 & 0 & \sigma v & -\sigma v & -\sigma v & \sigma v & 0 & 0 & 0 & 0 & 0 & 0\\
v & v & -v & -v & -v & -v & 0 & 0 & 0 & 0 & 0 & 0\\
0 & 0 & 0 & 0 & 0 & 0 & 0 & 0 & v & v & 0 & 0\\
0 & 0 & 0 & 0 & 0 & 0 & 0 & 0 & 0 & 0 & v & v
\end{smallmatrix}\right)\left(\begin{smallmatrix}z(1,t)\\ z(0,t)\end{smallmatrix}\right),\\
y(t) &= \left(\begin{smallmatrix}0 & 0 & 0 & 0 & \sigma v & -\sigma v & 0 & 0 & 0 & 0 & 0 & 0\\
0 & 0 & 0 & 0 & 0 & 0 & 0 & 0 & \sigma v & -\sigma v & 0 & 0
\end{smallmatrix}\right)\left(\begin{smallmatrix}z(1,t)\\ z(0,t)\end{smallmatrix}\right).
\end{align*}
In addition, let $\tilde{z}_i(\zeta,t) = z_i(1-\zeta,t), i = 1, 3, 5$. By denoting $\mathbf{z} := \left(\begin{smallmatrix}\tilde{z}_1 & z_2 & \tilde{z}_3 & z_4 & \tilde{z}_5 & z_6\end{smallmatrix}\right)^T$, there holds
\begin{align}
\frac{\partial\mathbf{z}}{\partial t}(\zeta,t) &= -v\frac{\partial\mathbf{z}}{\partial\zeta}(\zeta,t),\label{Dynamics_Strings}\\
\left(\begin{smallmatrix}0_{3\times 3}\\ I_{3\times 3}\end{smallmatrix}\right)u(t) &= -v\left(\begin{smallmatrix}0 & \sigma & 0 & 0 & 0 & 0\\
-\sigma & 0 & \sigma & 0 & 0 & 0\\
0 & 0 & -\sigma & 0 & \sigma & 0\\
-1 & 0 & 1 & 0 & 1 & 0\\
0 & 0 & 0 & -1 & 0 & 0\\
0 & 0 & 0 & 0 & 0 & -1
\end{smallmatrix}\right)\mathbf{z}(0,t) - v \left(\begin{smallmatrix}-\sigma & 0 & 0 & 0 & 0 & 0\\
0 & \sigma & 0 & -\sigma & 0 & 0\\
0 & 0 & 0 & \sigma & 0 & -\sigma\\
0 & -1 & 0 & 1 & 0 & 1\\
0 & 0 & -1 & 0 & 0 & 0\\
0 & 0 & 0 & 0 & -1 & 0
\end{smallmatrix}\right)\mathbf{z}(1,t)\nonumber\\
&= - v K\mathbf{z}(0,t) - v L\mathbf{z}(1,t),\label{Inputs_Strings}\\
y(t) &= -v\left(\begin{smallmatrix}0 & 0 & 0 & 0 & -\sigma & 0\\
0 & 0 & 0 & \sigma & 0 & 0
\end{smallmatrix}\right)\mathbf{z}(0,t) -v\left(\begin{smallmatrix}
0 & 0 & 0 & 0 & 0 & \sigma\\
0 & 0 & -\sigma & 0 & 0 & 0
\end{smallmatrix}\right)\mathbf{z}(1,t)\nonumber\\
&= - v K_y\mathbf{z}(0,t) - v L_y\mathbf{z}(1,t).\label{Outputs_Strings}
\end{align}
The determinant of $K$ is given by $\sigma^3$, which makes $K$ invertible if and only if $\sigma\neq 0$. This is always the case since $\sigma = (\rho T)^{-\frac{1}{2}}$. Hence, \eqref{Dynamics_Strings}--\eqref{Outputs_Strings} is a well-posed boundary controlled system, see Proposition \ref{Prop_Well_Posed}. According to the definition of the matrices $A_d, B_d, C_d$ and $D_d$ in \eqref{Matrices_Discrete}, there holds
\begin{align*}
A_d &= -K^{-1}L = \left(\begin{smallmatrix}0 & 3 & 0 & -2 & 0 & -2\\1 & 0 & 0 & 0 & 0 & 0\\
0 & 2 & 0 & -1 & 0 & -2\\
0 & 0 & -1 & 0 & 0 & 0\\
0 & 2 & 0 & -2 & 0 & -1\\
0 & 0 & 0 & 0 & -1 & 0\end{smallmatrix}\right),\\
B_d &= -K^{-1}\left[\begin{smallmatrix}0\\I\end{smallmatrix}\right] = \left(\begin{smallmatrix}-1 & 0 & 0\\
0 & 0 & 0\\
-1 & 0 & 0\\
0 & 1 & 0\\
-1 & 0 & 0\\
0 & 0 & 1\end{smallmatrix}\right),\\
C_d &= K_yK^{-1}L - L_y = \left(\begin{smallmatrix}0 & 2\sigma & 0 & -2\sigma & 0 & -2\sigma\\
0 & 0 & 2\sigma & 0 & 0 & 0\end{smallmatrix}\right),\\
D_d &= K_yK^{-1}\left[\begin{smallmatrix}0\\I\end{smallmatrix}\right] = \left(\begin{smallmatrix}-\sigma & 0 & 0\\
0 & -\sigma & 0\end{smallmatrix}\right).
\end{align*}
Note that the eigenvalues of the matrix $A_d$ are given by $i, -i, -1-\sqrt{2}, 1-\sqrt{2}, -1+\sqrt{2}, 1+\sqrt{2}$, which means that the open-loop system is not asymptotically stable. For simulation purposes, let us consider $v = 1$ and $\sigma = 1$. By using Matlab$^\copyright$\, with the matrices $K, L, K_y$ and $L_y$ defined in \eqref{Inputs_Strings} and \eqref{Outputs_Strings}, we get that the CARE \eqref{CARE} and the FARE \eqref{FARE} possess the following nonnegative self-adjoint solutions
\begin{align*}
\Pi &= \left(\begin{smallmatrix}1.6415 & 0 & 0 & 0 & 0.5702 & 0\\
0 & 5.6355 & 0 & -4.1131 & 0 & -3.5224\\
0 & 0 & 2 & 0 & 0 & 0\\
0 & -4.1131 & 0 & 3.8635 & 0 & 2.2497\\
0.5702 & 0 & 0 & 0 & 0.7067 & 0\\
0 & -3.5224 & 0 & 2.2497 & 0 & 3.2728
\end{smallmatrix}\right),\\
\tilde{\Pi} &= \left(\begin{smallmatrix}
1.6418 & 0 & 1.1988 & 0 & 0.9430 & 0\\
0 & 0.8465 & 0 & 0 & 0 & -0.6830\\
1.1988 & 0 & 1.3069 & 0 & 0.3919 & 0\\
0 & 0 & 0 & 0.5000 & 0 & 0\\
0.9430 & 0 & 0.3919 & 0 & 1.0511 & 0\\
0 & -0.6830 & 0 & 0 & 0 & 1.9661
\end{smallmatrix}\right),
\end{align*}
which makes the system \eqref{Dynamics_Strings}--\eqref{Outputs_Strings} input and output stabilizable, see Remark \ref{Charact_IO_Stabilizab}. Observe now that the operator $F_d$ defined in \eqref{Optimal_Input_Original_Var} is given by
\begin{align*}
F_d = \left(\begin{smallmatrix}0 & 2.0283 & 0 & -1.4659 & 0 & -1.5624\\
0.4827 & 0 & 1 & 0 & 0.1125 & 0\\
0.5702 & 0 & 0 & 0 & 0.7067 & 0
\end{smallmatrix}\right) =: \left(\begin{smallmatrix}0 & \beta_1 & 0 & \beta_2 & 0 & \beta_3\\
\beta_4 & 0 & 1 & 0 & \beta_5 & 0\\
\beta_6 & 0 & 0 & 0 & \beta_7 & 0
\end{smallmatrix}\right),
\end{align*}
which implies that $\sigma(A_d+B_dF_d)=\{0.6839, -0.6839, 0.4780i, -0.4780i, 0, 0\}$. Hence, the matrix $A_d+B_dF_d$ is strongly stable. Thanks to Theorem \ref{Thm_Unique_Sol}, $\Pi$ is the unique nonnegative, self-adjoint and stabilizing solution of the CARE \eqref{CARE}. According to \eqref{Optimal_Input_Original_Var}, the optimal feedback operator is given by 
\begin{align*}
u^\text{opt}(t) &= F_d\mathbf{z}^\text{opt}(1,t)\\
			 &= \left(\begin{smallmatrix}0 & \beta_1 & 0 & \beta_2 & 0 & \beta_3\\
\beta_4 & 0 & 1 & 0 & \beta_5 & 0\\
\beta_6 & 0 & 0 & 0 & \beta_7 & 0
\end{smallmatrix}\right)\left(\begin{smallmatrix}z_1^\text{opt}(0,t)\\ z_2^\text{opt}(1,t)\\ z_3^\text{opt}(0,t)\\ z_4^\text{opt}(1,t)\\ z_5^\text{opt}(0,t)\\ z_6^\text{opt}(1,t)\end{smallmatrix}\right)\\
&= \frac{1}{2}\left(\begin{smallmatrix}0 & \beta_1 & 0 & \beta_2 & 0 & \beta_3\\
\beta_4 & 0 & 1 & 0 & \beta_5 & 0\\
\beta_6 & 0 & 0 & 0 & \beta_7 & 0
\end{smallmatrix}\right)\left(\begin{smallmatrix}x_1^\text{opt}(0,t) + \sigma^{-1} x_2^\text{opt}(0,t)\\ \sigma^{-1} x_2^\text{opt}(1,t)-x_1^\text{opt}(1,t)\\ x_3^\text{opt}(0,t)+\sigma^{-1} x_4^\text{opt}(0,t)\\ \sigma^{-1} x_4^\text{opt}(1,t)-x_3^\text{opt}(1,t)\\ x_5^\text{opt}(0,t)+\sigma^{-1} x_6^\text{opt}(0,t)\\ \sigma^{-1} x_6^\text{opt}(1,t)-x_5^\text{opt}(1,t)\end{smallmatrix}\right).
\end{align*}
In terms of the original state variables, the optimal control inputs $u^\text{opt}_1(t), u^\text{opt}_2(t)$ and $u^\text{opt}_3(t)$ are given by
\begin{align*}
&u^\text{opt}_1(t) =\\
&\frac{1}{2}\left[\beta_1\left(\sigma^{-1}\frac{\partial w_\text{A}^\text{opt}}{\partial\zeta}(1,t)-\rho\frac{\partial w^\text{opt}_\text{A}}{\partial t}(1,t)\right) + \beta_2\left(\sigma^{-1}\frac{\partial w_\text{B}^\text{opt}}{\partial\zeta}(1,t)-\rho\frac{\partial w^\text{opt}_\text{B}}{\partial t}(1,t)\right)\right.\\
&\left. + \beta_3\left(\sigma^{-1}\frac{\partial w^\text{opt}_\text{C}}{\partial\zeta}(1,t)-\rho\frac{\partial w^\text{opt}_\text{C}}{\partial t}(1,t)\right)\right],
\end{align*}
\begin{align*}
&u^\text{opt}_2(t) =\\
&\frac{1}{2}\left[\beta_4\left(\sigma^{-1}\frac{\partial w_\text{A}^\text{opt}}{\partial\zeta}(0,t)+\rho\frac{\partial w^\text{opt}_\text{A}}{\partial t}(0,t)\right) + \left(\sigma^{-1}\frac{\partial w^\text{opt}_\text{B}}{\partial\zeta}(0,t)+\rho\frac{\partial w^\text{opt}_\text{B}}{\partial t}(0,t)\right)\right.\\
&\left. + \beta_5\left(\sigma^{-1}\frac{\partial w^\text{opt}_\text{C}}{\partial\zeta}(0,t)+\rho\frac{\partial w^\text{opt}_\text{C}}{\partial t}(0,t)\right)\right],
\end{align*}
and
\begin{align*}
&u^\text{opt}_3(t)=\\
&\frac{1}{2}\left[\beta_6\left(\sigma^{-1}\frac{\partial w_\text{A}^\text{opt}}{\partial\zeta}(0,t)+\rho\frac{\partial w_\text{A}^\text{opt}}{\partial t}(0,t)\right) + \beta_7\left(\sigma^{-1}\frac{\partial w^\text{opt}_\text{C}}{\partial\zeta}(0,t)+\rho\frac{\partial w^\text{opt}_\text{C}}{\partial t}(0,t)\right)\right],
\end{align*}
respectively.

\subsection{Co-current heat exchanger}
\begin{figure*}
\begin{center}
%\documentclass[border=1pt,tikz]{standalone}
%\usepackage{amsmath} % for \dfrac
%\usepackage{bm}
%\usepackage{physics}
%\usepackage{tikz,pgfplots}
%\usetikzlibrary{angles,quotes} % for pic (angle labels)
%\usetikzlibrary{calc}
%\usetikzlibrary{decorations.markings}
%\tikzset{>=latex} % for LaTeX arrow head
%
%\usepackage{xcolor}
\colorlet{Ecol}{orange!90!black}
\colorlet{EcolFL}{orange!80!black}
\colorlet{veccol}{green!45!black}
\colorlet{EFcol}{red!60!black}
\tikzstyle{charged}=[top color=blue!20,bottom color=blue!40,shading angle=10]
\tikzstyle{darkcharged}=[very thin,top color=blue!60,bottom color=blue!80,shading angle=10]
\tikzstyle{charge+}=[very thin,top color=red!80,bottom color=red!80!black,shading angle=-5]
\tikzstyle{charge-}=[very thin,top color=blue!50,bottom color=blue!70!white!90!black,shading angle=10]
\tikzstyle{gauss surf}=[red!40!black,top color=green!2,bottom color=red!80!black!70,shading angle=5,fill opacity=0.5]
\tikzstyle{gauss lid}=[gauss surf,middle color=red!80!black!20,shading angle=40,fill opacity=0.6]
\tikzstyle{gauss dark}=[red!50!black,fill=red!60!black!70,fill opacity=0.8]
\tikzstyle{gauss line}=[red!40!black]
\tikzstyle{gauss dashed line}=[red!60!black!80,dashed,line width=0.1]
\tikzstyle{EField}=[->,thick,Ecol]
\tikzstyle{vector}=[->,thick,veccol]
\tikzstyle{normalvec}=[->,thick,blue!80!black!80]
\tikzstyle{EFieldLine}=[thick,EcolFL,decoration={markings,
          mark=at position 0.5 with {\arrow{latex}}},
          postaction={decorate}]
\tikzstyle{measure}=[fill=white,midway,outer sep=2]
\def\L{2}
\def\W{0.2}
\def\N{4}
%
%
%
%\begin{document}
\begin{tikzpicture}[scale=4]
  \def\M{8}
  \def\R{0.4*\L}
  \def\g{0.2*\R}
  \def\G{0.4*\R}
  \def\a{0.33*\L}
  \coordinate (L)  at (-\a,0);
  \coordinate (R)  at (+\a,0);
  \coordinate (TL) at (-\a,\G);
  \coordinate (TR) at (+\a,\G);
  \coordinate (BL) at (-\a,-\G);
  \coordinate (BR) at (+\a,-\G);
  
  % GAUSS BEHIND
  \draw[gauss line] (TR) arc (90:270:{\g} and {\G});
  
  % ROD
  \draw[charged] (-6.35*\L/16,-\W/2) --++(5.8*\L/8,0) to[out=0,in=0] ++ (0,\W) --++ (-5.8*\L/8,0) -- cycle;
  \draw[charged] (-6.35*\L/16,-\W/2) to[out=180,in=180] ++ (0,\W) to[out=0,in=0] cycle;
  \draw[->] (-5*\L/16,0) -- (-2.5*\L/16,0);
  \node[] at (0.07,0) {$T_i(\zeta,t)$};
  \draw[->] (3.9*\L/16,0) -- (6.4*\L/16,0);
  
  \begin{scope}
    \clip (-\L/2,-0.5*\W)
      --++ (\L/2-\a,0) to[out=50,in=-50] ++(0,1.0*\W) --++ (-\L/2+\a,0) --++ (0,\G)
      --++ (\L-2*\a,0) --++ (0,{-2*(\G+\W)}) --++ (-\L+2*\a,0) -- cycle;
    \draw[gauss lid] (L) ellipse ({\g} and {\G});
  \end{scope}

  % GAUSS IN FRONT
  %\draw[<->] (-\a,\W/2) -- (-\a,\G) node[measure,fill=green!80!black!8,inner sep=1,outer sep=0] {$r$};
  \draw[|->|,dashed] (-\a,-1.5*\G) -- (\a,-1.5*\G);% node[measure] {};
  \node[above] at (0,-1.5*\G) {$\zeta$};
  \node[below] at (-\a,-1.5*\G) {$0$};
  \node[below] at (\a,-1.5*\G) {$1$};
  \draw[->] (-\a,1.5*\G) -- (-\a,1*\G);
  \node[above] at (-\a,1.6*\G) {$u(t)$};
  %\draw[|-|] (-\a,0) -- (\a,0) node[measure] {\tiny{$T_i(\zeta,t)$}};
  %\draw[normalvec] (-1.1*\a,-0.6*\G) --++ (-0.25*\G,0) node[left] {$\vu{n}$};
  %\draw[normalvec] (+1.1*\a,-0.6*\G) --++ (+0.25*\G,0) node[above] {$\vu{n}$};
  \draw[gauss surf]
    (BL) arc (-90:90:{\g} and {\G}) --
    (TR) arc (90:-90:{\g} and {\G}) -- cycle;
   \draw[->] (-5*\L/16,0.2) -- (-2.5*\L/16,0.2);
  \node[] at (0.07,0.2) {$T_e(\zeta,t)$};
  \draw[->] (3.9*\L/16,0.2) -- (6*\L/16,0.2);
  %\draw[normalvec] (-0.064*\L,\G) --++ (0,0.25*\G) node[left=2,above] {$\vu{n}$};
  %\draw[->,thick,Ecol] (-0.045*\L,\G) --++ (0,0.6*\G) node[right] {$\vb{E}$};
  
  % LABELS
  %\node[left,green!30!black] at ($(-\a,0)+(150:{\g} and {\G})$) {$S_1$};
  %\node[right,green!30!black] at ($(\a,0)+(30:{\g} and {\G})$) {$S_2$};
  %\node[above,green!30!black] at (0.6*\a,\G) {$S_3$};
  
\end{tikzpicture}
%\end{document}
\end{center}
\caption{Schematic profile view of a co-current heat-exchanger.\label{fig:HE}}
\end{figure*}
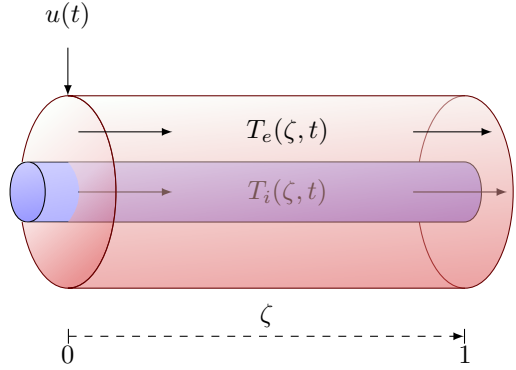

In this part, we consider a system composed of a co-current heat exchanger of length $1$, as depicted in Figure \ref{fig:HE}.

The device consists in two tubes, one immersed in the other. A fluid flows in each tube and the heat exchange occurs only between the two tubes, meaning that we neglect the heat losses with the external environment. The temperature in the internal and the external tubes at position $\zeta\in [0,1]$ and time $t\geq 0$ are denoted by $T_i(\zeta,t)$ and $T_e(\zeta,t)$, respectively. By defining $\mathbf{T}(\zeta,t) = (\begin{smallmatrix}T_i(\zeta,t) & T_e(\zeta,t)\end{smallmatrix})^T$, the PDE governing this system is given by
\begin{equation}
\frac{\partial \mathbf{T}}{\partial t}(\zeta,t) = -v(\zeta)\frac{\partial\mathbf{T}}{\partial \zeta}(\zeta,t) + \left(\begin{matrix}-\alpha_1(\zeta) & \alpha_1(\zeta)\\ \alpha_2(\zeta) & -\alpha_2(\zeta)\end{matrix}\right)\mathbf{T}(\zeta,t),
\label{PDE_Temp}
\end{equation}
where $\alpha_i(\zeta) > 0, i=1, 2, \forall \zeta\in [0,1]$, are heat transfer functions and $v(\zeta) > 0, \forall \zeta\in [0,1]$, is the velocity at which the fluid is flowing, at position $\zeta$. In what follows, we denote by $M(\zeta)$ the matrix $ \left(\begin{smallmatrix}-\alpha_1(\zeta) & \alpha_1(\zeta)\\ \alpha_2(\zeta) & -\alpha_2(\zeta)\end{smallmatrix}\right)$. The initial temperature profile is given by $\mathbf{T}(\zeta,0) = \mathbf{T}_0(\zeta)$. Such a setting may be found for instance in \cite{Maidi_2010}. To the PDE \eqref{PDE_Temp}, we add the following boundary condition
\begin{equation}
T_i(0,t) = 0, t\geq 0.
\label{BC_HE}
\end{equation}
The temperature at the inlet of the external tube is taken as control input and the measured output is the difference temperature between the outlet and the inlet of the external tube. In other words, 
\begin{align}
u(t) &= T_e(0,t)\label{Input_HE},\\
y(t) &= T_e(1,t) - T_e(0,t)\label{Output_HE},
\end{align}
$t\geq 0$. Hence, \eqref{PDE_Temp} together with \eqref{BC_HE}--\eqref{Output_HE} may be written as
\begin{align}
\frac{\partial \mathbf{T}}{\partial t}(\zeta,t) &= -v(\zeta)\frac{\partial\mathbf{T}}{\partial \zeta}(\zeta,t) + M(\zeta)\mathbf{T}(\zeta,t),\label{State_Co_Current}\\
\left(\begin{matrix}0\\ 1\end{matrix}\right)u(t) &= \mathbf{T}(0,t),\label{Input_Co_Current}\\
y(t) &= (\begin{matrix}0 & -1\end{matrix})\mathbf{T}(0,t) + (\begin{matrix}0 & 1\end{matrix})\mathbf{T}(1,t).\label{Output_Co_Current}
\end{align}
In order for \eqref{State_Co_Current}--\eqref{Output_Co_Current} to be expressed in the same way as \eqref{StateEquation_Diag_Uniform}--\eqref{Output_Diag_Uniform}, we define the new state variable $\mathcal{T}(\zeta,t)$ as $\mathcal{T}(\zeta,t) := \frac{1}{v(\zeta)}\mathbf{T}(\zeta,t)$. By dividing the PDE \eqref{State_Co_Current} by $v(\zeta)$ and by re-expressing \eqref{Input_Co_Current}--\eqref{Output_Co_Current} with $\mathcal{T}(\zeta,t)$, we get that 
\begin{align}
\frac{\partial \mathcal{T}}{\partial t}(\zeta,t) &= -\frac{\partial}{\partial \zeta}(v(\zeta)\mathcal{T}(\zeta,t)) + M(\zeta)\mathcal{T}(\zeta,t),\label{State_Co_Current_Def}\\
\left(\begin{matrix}0\\ 1\end{matrix}\right)u(t) &= v(0)\mathcal{T}(0,t),\label{Input_Co_Current_Def}\\
y(t) &= (\begin{matrix}0 & -1\end{matrix})v(0)\mathcal{T}(0,t) + (\begin{matrix}0 & 1\end{matrix})v(1)\mathcal{T}(1,t).\label{Output_Co_Current_Def}
\end{align}
Hence it is of the form \eqref{StateEquation_Diag_Uniform}--\eqref{Output_Diag_Uniform} with
\begin{align*}
\lambda_0(\cdot) = v(\cdot),\,\,\, K = -I,\,\,\, L = 0,\,\,\, K_y = (\begin{matrix}0 & 1\end{matrix}),\,\,\, L_y = -(\begin{matrix}0 & 1\end{matrix}).
\end{align*}
As $K$ is invertible, \eqref{State_Co_Current_Def}--\eqref{Output_Co_Current_Def} is a well-posed linear boundary controlled system. In order to be able to compute the matrices $A_d, B_d, C_d$ and $D_d$ in \eqref{Matrices_Discrete}, let us have a look at the differential equation 
\begin{align*}
P'(\zeta) = \frac{1}{v(\zeta)}M(\zeta)P(\zeta),\hspace{0.5cm} P(0) = I,
\end{align*}
see \eqref{Equa_Diff_P}. The solution of this differential equation is given by
\begin{align*}
P(\zeta) = \left(\begin{matrix}\displaystyle 1-\int_0^\zeta e^{h(\eta)}\frac{\alpha_1(\eta)}{v(\eta)}\d\eta & \displaystyle 1 - e^{h(\zeta)} - \int_0^\zeta e^{h(\eta)}\frac{\alpha_2(\eta)}{v(\eta)}\d\eta\\
\displaystyle 1 - e^{h(\zeta)}-\int_0^\zeta e^{h(\eta)}\frac{\alpha_1(\eta)}{v(\eta)}\d\eta & \displaystyle 1 - \int_0^\zeta e^{h(\eta)}\frac{\alpha_2(\eta)}{v(\eta)}\d\eta\end{matrix}\right),
\end{align*}
with
\begin{align*}
\displaystyle h:[0,1]\to \mathbb{R}^-, h(\eta) := -\int_0^\eta\left[\frac{\alpha_1(\tau)+\alpha_2(\tau)}{v(\tau)}\right]\d\tau.
\end{align*}
As a consequence, the matrices $A_d, B_d, C_d$ and $D_d$ in \eqref{Matrices_Discrete} are given by
\begin{align*}
A_d &= -K^{-1}LP(1) = 0,\\
B_d &= -K^{-1}\left(\begin{matrix}0\\ 1\end{matrix}\right) = \left(\begin{matrix}0\\ 1\end{matrix}\right),\\
C_d &= (K_yK^{-1}L-L_y)P(1)\\
&= (\begin{matrix}0 & 1\end{matrix})P(1) = \left(\begin{smallmatrix}\displaystyle 1-e^{h(1)} -\int_0^1 e^{h(\eta)}\frac{\alpha_1(\eta)}{v(\eta)}\d\eta &\hspace{0.7cm} \displaystyle 1 - \int_0^1 e^{h(\eta)}\frac{\alpha_2(\eta)}{v(\eta)}\d\eta\end{smallmatrix}\right),\\
D_d &= K_yK^{-1}\left(\begin{matrix}0\\ 1\end{matrix}\right) = -1.
\end{align*}
As the matrix $A_d = 0$, the discrete-time system $(A_d,B_d,C_d,D_d)$ is input and output stable, and hence input and output stabilizable, which makes the finite cost conditions satisfied, see \cite[Lemmas 3.3 and 3.6]{Opmeer_LQ_Discrete}. As a consequence, it is also known that the control algebraic Riccati equation \eqref{CARE} has a nonnegative self-adjoint solution. In our setting, the CARE \eqref{CARE} takes the following form
\begin{equation*}
-\Pi + C_d^*C_d = \left(1+D_d^*D_d+B_d^*\Pi B_d\right)^{-1}C_d^*C_d,
\end{equation*}
which can be equivalently written as
\begin{equation*}
\Pi = C_d^*C_d\left(1-\left(2+\left(\begin{smallmatrix}0 & 1\end{smallmatrix}\right)\Pi\left(\begin{smallmatrix}0 \\ 1\end{smallmatrix}\right)\right)^{-1}\right).
\end{equation*}
Its nonnegative solution is given by
\begin{align}
\Pi &= \left(\begin{matrix}\pi_1 & \tilde{\pi}\\ \tilde{\pi} & \pi_2\end{matrix}\right),\label{Sol_CARE_HE}
\end{align}
in which 
\begin{align*}
\pi_2 &= \frac{r_2 - 2 + \sqrt{r_2^2+4}}{2},\\
\pi_1 &= r_1(1-(2+\pi_2)^{-1}),\\
\tilde{\pi} &= \tilde{r}(1-(2+\pi_2)^{-1}),
\end{align*}
with
\begin{align*}
\left(\begin{matrix}r_1 & \tilde{r}\\\tilde{r} & r_2\end{matrix}\right) := C_d^*C_d,
\end{align*}
that is,
\begin{align*}
r_1 &= \displaystyle\left(1-e^{h(1)}-\int_0^1 e^{h(\eta)}\frac{\alpha_1(\eta)}{v(\eta)}\d\eta\right)^2,\\
\tilde{r} &= \displaystyle\left(1-e^{h(1)}-\int_0^1 e^{h(\eta)}\frac{\alpha_1(\eta)}{v(\eta)}\d\eta\right)\left(1 - \int_0^1 e^{h(\eta)}\frac{\alpha_2(\eta)}{v(\eta)}\d\eta\right),\\
r_2 &= \displaystyle\left(1 - \int_0^1 e^{h(\eta)}\frac{\alpha_2(\eta)}{v(\eta)}\d\eta\right)^2.
\end{align*}
Now we compute the optimal feedback operator $F_d$. According to \eqref{Optimal_Input_Original_Var}, it is given by
\begin{align*}
F_d &= -(I + D_d^*D_d + B_d^*\Pi B_d)^{-1}(D_d^*C_d + B^*_d\Pi A_d)\\
&= (2+\pi_2)^{-1}\left(\begin{smallmatrix}\displaystyle 1-e^{h(1)} -\int_0^1 e^{h(\eta)}\frac{\alpha_1(\eta)}{v(\eta)}\d\eta &\hspace{0.7cm} \displaystyle 1 - \int_0^1 e^{h(\eta)}\frac{\alpha_2(\eta)}{v(\eta)}\d\eta\end{smallmatrix}\right).
\end{align*}
Hence, the closed-loop operator dynamics $A_\Pi$ is expressed by 
\begin{equation*}
A_\Pi = B_dF_d = (2+\pi_2)^{-1}\left(\begin{matrix}0 & 0\\ \mu & \kappa\end{matrix}\right),
\end{equation*}
where $\mu = 1-e^{h(1)} -\int_0^1 e^{h(\eta)}\frac{\alpha_1(\eta)}{v(\eta)}\d\eta$ and $\kappa$ is such that $\kappa^2 = r_2$.
It is easy to see that the eigenvalues of the matrix $A_\Pi$ are given by 
\begin{align*}
\lambda_1 = 0, \lambda_2 = \frac{\kappa}{2+\pi_2}.
\end{align*}
According to the expression of $\pi_2$, there holds $\vert\lambda_2\vert < 1$, which implies that the matrix $A_\Pi$ is stable in the sense that its spectral radius is less than 1. As a consequence, the solution $\Pi$ \eqref{Sol_CARE_HE} to the CARE \eqref{CARE} is unique. The optimal control input in continuous-time, $u^\text{opt}(t)$, is given by
\begin{align*}
u^\text{opt}(t) &= \lambda_0(1)(2+\pi_2)^{-1}C_dQ(1)\mathcal{T}^\text{opt}(1,t)\\
&= (2+\pi_2)^{-1}(\begin{matrix}0 & 1\end{matrix})P(1)Q(1)\mathbf{T}^\text{opt}(1,t)\\
&= (2+\pi_2)^{-1}T_e^{\text{opt}}(1,t).
\end{align*}

\section{Perspectives and future works}
\label{sec:ccl}
As a first perspective, the extension of the proposed approach to systems driven by the PDEs of the form
\begin{equation}
\frac{\partial z}{\partial t}(\zeta,t) = \frac{\partial}{\partial\zeta}\left(\left(\begin{smallmatrix}\lambda(\zeta) & 0\\ 0 & -\lambda(\zeta)\end{smallmatrix}\right)z(\zeta,t)\right) + M(\zeta)z(\zeta,t),
\label{ExtendedSystem}
\end{equation}
is worth investigating. Indeed, the class of systems \eqref{ExtendedSystem} pops up for instance when looking at the dynamics of a Timoshenko beam after writing the system in Riemann coordinates. The class \eqref{ExtendedSystem} also appears when describing the dynamics of counter-current heat exchangers. The main challenge with \eqref{ExtendedSystem} is that the matrix $M$ may not be removed by a change of variables similar to the one introduced in Lemma \ref{Lemma_Chg_Var}. This would be the case provided that the matrices $\left(\begin{smallmatrix}\lambda & 0\\ 0 & -\lambda\end{smallmatrix}\right)$ and $M$ commute, which holds under strong assumptions like the diagonal property of $M$ for instance. Trying to extend our approach to the case where the function $\lambda$ is replaced by a diagonal matrix with different entries is also a challenge for future works. This would greatly enlarge the applicability of our results. As is extensively shown in the book \cite{BastinCoron_Book}, the PDEs \eqref{StateEquation_Diag_Uniform} with $\lambda$ replaced by a diagonal matrix with different entries are able to model many physical phenomena. One particular example is the movement of water and sediment in open channels, governed by the Saint-Venant-Exner equations. Extending our approach to the class of port-Hamiltonian systems considered in \cite{JacobZwart} would be interesting to look at. This class of systems is able to model many physical phenomena and it possesses a lot of interesting and useful properties. To mention a few, the interconnection of two port-Hamiltonian systems is again a port-Hamiltonian system. Moreover, the energy in the system is related to the $\L^2$-norm of the state trajectory, which is useful for looking at system properties such as stability for instance. Those properties could be exploited to solve the LQ-optimal control problem for this class of systems. As an additional perspective, numerical simulations could be depicted to highlight the performances of the proposed method. In addition, future works would aim at investigating the dual approach to LQ-optimal control in order to do optimal state reconstruction, that is, Kalman filtering. The meaning of duality in this context should be properly addressed because of the unboundedness of the control and the observation operators. LQ-optimal control can be combined with a state estimator leading to the deterministic version of the Kalman filter, see \cite[Exercises 9.23, 9.24]{CurtainZwart2020}. The fact that our system may be transformed into a discrete-time system should lead to simplifications of this problem.

\bibliographystyle{siamplain}
\bibliography{biblio_Main}

\end{document}